\documentclass{article}
\usepackage{amsmath}
\usepackage{amssymb}
\usepackage{graphics,graphicx,psfrag}

\newenvironment{proof}{\noindent {\bf Proof}}
{\hfill $\square$ \vspace{0.5cm}}

\newcommand{\vip}{\vskip0.15cm}

\newcommand{\nn}{\mathbb{N}}
\newcommand{\zz}{\mathbb{Z}}
\newcommand{\rr}{{\mathbb{R}}}

\newcommand{\epr}{{\{0,1\}^\zz}}

\newcommand{\cL}{{\cal L}}
\newcommand{\cP}{{\cal P}}
\newcommand{\cG}{{\cal G}}
\newcommand{\cH}{{\cal H}}
\newcommand{\cA}{{\cal A}}

\newcommand{\cF}{{\cal F}}

\newcommand{\tzeta}{{\tilde \zeta}}
\newcommand{\bzeta}{{\bar \zeta}}
\newcommand{\bareta}{{\bar \eta}}

\newcommand{\tZ}{{\tilde Z}}

\newcommand{\tR}{{\tilde R}}
\newcommand{\tT}{{\tilde T}}

\newcommand{\oR}{{\overline R}}
\newcommand{\uL}{{\underline L}}

\newcommand{\bZ}{{\mathbf 0}}

\newcommand{\sm}{{{s-}}}

\newcommand{\ala}{\nonumber \\}
\newcommand{\indiq}{1\!\! 1}
\newcommand{\e}{{\varepsilon}}

\newcommand{\hzeta}{\widehat{\zeta}}
\newcommand{\heta}{\widehat{\eta}}

\newtheorem{theo}{\indent Theorem}[section]
\newtheorem{prop}[theo]{\indent Proposition}

\newtheorem{lem}[theo]{\indent Lemma}
\newtheorem{defin}[theo]{\indent Definition}
\newtheorem{cor}[theo]{\indent Corollary}
\newtheorem{nota}[theo]{\indent Notation}

\newtheorem{gc}[theo]{\indent Graphical construction}

\begin{document}

\title{On the invariant distribution of a one-dimensional avalanche process}

\author{Xavier {\sc Bressaud}\footnote{IML, Case 907, 163 avenue de Luminy, 
13288 Marseille Cedex 9, France. E-mail {\tt bressaud@iml.univ-mrs.fr}}\;
and Nicolas {\sc Fournier}\footnote{Centre de Maths,
Facult\'e des Sciences et Technologies,
Universit\'e Paris~12, 61 avenue du G\'en\'eral de Gaulle, 94010 Cr\'eteil 
Cedex, France. E-mail: {\tt nicolas.fournier@univ-paris12.fr} }}

\maketitle

\def\abstractname{Abstract}
\begin{abstract}
\noindent We consider an interacting particle system $(\eta_t)_{t\geq 0}$
with values in $\epr$,
in which each vacant site becomes occupied with rate $1$, while each 
connected component of occupied sites become vacant with  rate equal to
its size. We show that such a process 
admits a unique invariant distribution, 
which is exponentially mixing and can be perfectly simulated.
We also prove that for any initial condition, the avalanche process
tends to equilibrium exponentially fast, as time increases to infinity.
Finally, we consider a related mean-field coagulation-fragmentation 
model, we compute its 
invariant distribution, and we show numerically that it is very close
to that of the interacting particle system.
\end{abstract}


\vip


\section{Notations and main results}
\label{sex1} \setcounter{equation}{0}

Consider an independent family $N=((N_t(i))_{t\geq 0})_{i \in \zz}$ 
of Poisson processes with rate $1$. In the whole paper, such
a family will be called an IFPP.

Assume that on each site $i\in\zz$, snow flocks are falling according to 
the process $(N_t(i))_{t\geq 0}$. When a flock falls on a vacant 
site of $\zz$, this site becomes
occupied. When a flock falls on an occupied site $i\in \zz$,
an {\it avalanche} starts:
the whole connected component of occupied sites around $i$ becomes vacant.

We denote by $(\eta_t(i))_{t\geq 0, i \in \zz}$ the process
defined,
for $t\geq 0$ and $i\in \zz$, by $\eta_t(i)=1$ 
(resp. $\eta_t(i)=0$) if the site $i$ is occupied (resp. vacant) at time $t$.

\vip

To avoid infinite rates of interaction, we will restrict our study
to the case where the initial condition $\eta_0$ lies in the following
space:
\begin{equation}
E:= \left\{\eta \in \epr; \;\; 
\liminf_{i \to -\infty} \eta(i)=\liminf_{i \to +\infty} \eta(i)=0  \right\}.
\end{equation}
A state $\eta$ belongs to $E$ if and only if it 
has no infinite connected component of occupied sites.
This condition is not really restrictive: easy considerations
show that even if $\eta_0 \in \epr \backslash E$, $\eta_t \in E$ for all 
$t>0$. This comes from the fact that infinite
connected components of occupied sites have an infinite death rate.

\vip

It is standard and easy to show, using for example a 
{\it graphical construction},
that for any initial condition $\eta_0 \in E$, for any IFPP $N$, 
the process $(\eta_t)_{t\geq 0}$ exists, is unique, and takes its values
in $E$.
It is actually a deterministic function of  $\eta_0$ and $(N(i))_{i\in\zz}$.
We call this process the $\eta_0$-avalanche process,
or the $(\eta_0,N)$-avalanche process when this precision is needed.
See \cite{liggett} for many examples of graphical constructions.

\vip

Furthermore, the process $(\eta_t)_{t\geq 0}$ is a strong Markov process,
and its infinitesimal generator $\cA$ is defined,
for $\eta \in E$ and $\Phi: E \mapsto \rr$ sufficiently regular 
(e.g. $\Phi(\eta)$ depending only on a finite number 
of coordinates of $\eta$) by
\begin{equation}
\cA \Phi (\eta) = \sum_{i\in\zz} \left[\Phi(a_i(\eta)) - \Phi(\eta)\right],
\end{equation}
where $a_i(\eta)\in E $ is defined in the following way:
\begin{itemize}
\item if $\eta(i)=0$, then $[a_i(\eta)](i)=1$ and  
$[a_i(\eta)](k)=\eta(k)$ for all 
$k \ne i$;
\item if $\eta(i)=1$, set $l_\eta(i)=\sup\{k \leq i; \eta(k)=0\}+1$,
$r_\eta(i)=\inf\{k \geq i; \eta(k)=0\}-1$, and put
$[a_i(\eta)](k)=0$ for $k\in [l_\eta(i),r_\eta(i)]$
and  $[a_i(\eta)](k)=\eta(k)$  for all $k \notin [l_\eta(i),r_\eta(i)]$.
\end{itemize}

Our main result in this paper concerns the invariant distribution
of the avalanche process.

For $A \subset \zz$ and
$\Gamma \in \cP(E)$
we denote by $\Gamma_A=\Gamma \circ p_A^{-1} \in \cP(\{0,1\}^A)$ 
its restriction to $A$, where $p_A: E \mapsto \{0,1\}^A$ is
the canonical projection.

For two probability measures $\mu$, $\nu$ 
on a measurable space $(F, \cF)$, 
we denote by
$|\mu - \nu|_{TV}= \sup_{G \in \cF} |\mu(G) - \nu(G)|$
the total variation between $\mu$ and $\nu$.

\begin{theo}\label{main}
(a) The avalanche process admits an unique invariant distribution 
$\Pi \in \cP(E)$.

(b) The exponential trend to equilibrium holds in the following sense. 
For $\varphi\in E$, denote by $\Pi_t^{\varphi}$ the
law of the $\varphi$-avalanche process at time $t$.
There exist some constants $C>0$, $\alpha>0$ such that
for all $t\geq 0$, all $l \geq 0$,
\begin{equation}\label{tteap}
\sup_{\varphi \in E} \left|(\Pi_t^{\varphi})_{[-l,l]} - \Pi_{[-l,l]}
\right|_{TV}
\leq C(1+l) e^{-\alpha t}.
\end{equation}
(c) For $l\geq 0$, 
there exists an explicit (see Appendix \ref{appendix}) 
and perfect simulation algorithm
for a $\Pi_{[-l,l]}$-distributed random variable $(\eta_0(i))_{i\in[-l,l]}$.

(d) The invariant distribution $\Pi$ is exponentially mixing 
in the following sense: one may find some constants $C>0$, $0<q<1$,
such that for any $k\in \zz$, $n \in \zz_+$,
\begin{equation}\label{mix}
\left|\Pi_{(-\infty,k]\cup[k+n,\infty)} - 
\Pi_{(-\infty,k]} \otimes \Pi_{[k+n,\infty)}  \right|_{TV} \leq C q^n.
\end{equation}
\end{theo}

Let us comment on these results. First, the system is very stable,
in the sense that no large clusters of occupied sites may appear. Indeed,
large clusters have a large death rate.
Clearly, the existence of invariant distributions should easily follow
from such an argument.
Of course, uniqueness of the invariant distribution and trend to equilibrium
are not suprising, but much more work is required, 
especially to give a rate of convergence.
The perfect simulation algorithm we give is quite complicated,
but gives, in some sense, an explicit expression
of the invariant distribution $\Pi$.
Finally, point (d) explains that at equilibrium, for two sites
$i$ and $j$, the dependance between $\eta(i)$ and $\eta(j)$ 
decreases exponentially fast with $|i-j|$. Such a result is 
also quite natural, but the proof is quite complicated.

\vip

At the end of the 80's, the so-called self-organized critical (SOC) systems
became rather popular. They are simple models supposed to enlight
temporal and spatial randomness observed in a variety of natural 
phenomena showing {\it long range correlations}, like  sand piles, avalanches,
earthquakes, stock market crashes, fire forest, shape of mountains,
of clouds, ... Very roughly,
the key idea (present in Bak-Tang-Wiesenfeld \cite{bak1} 
about sand piles) is that of
systems {\it growing} toward a {\it critical state}  
and relaxing through {\it catastrophic}  
events (avalanches, crashes, fire, ...)~; if the catastrophic 
events become more and more probable when approaching the critical state, 
the system spontaneously reaches an equilibrium {\it close} to the critical 
state.  

SOC systems commonly share other features as long range correlations, 
power laws for the amplitude of {\it catastrophic} events, 
spatial fractality of observed patterns, lack of typical scale, ... 
The most classical model is the so-called sand pile model introduced 
in 1987 in \cite{bak1}, but a lot of variants 
or related models have been proposed and studied more or less rigorously,  
describing earthquakes (Olami-Feder-Christensen, \cite{olami}) or 
fire forest (Henley \cite{henley}~; Drossel-Schwabl, \cite{drossel}) 
to mention a few. 
For surveys on the subject, see \cite{dhar} or \cite{bak2}, for instance.  

\vip

Initially, 
our process was thought as a very rough simplification of a sand pile model.  
In sand piles geometric rules describe the structure of a {\it stable} 
sand pile.  
Sand grains fall on a given pile~; if the new pile is {\it unstable},  it is 
re-organized to become stable, through (possibly many successive) 
elementary steps~; such events are called {\it avalanches}.  If the 
pile lives on a bounded domain, grains falling out of the domain 
disappear~; if the model is realistic, one can imagine that the 
number of grains in the pile and the shape of the pile reaches 
an equilibrium. Frequence and amplitude of {\it avalanches} at 
equilibrium are related to the number of grains that disappear. 
In our much simpler model, a grain falling on an occupied
site yields an {\it avalanche} involving all grains in the connected
component (that immediately disappear). It does not pretend to be a 
good physical description of a sand pile: the purpose is more to 
catch what is really  important in SOC systems. 

\vip

This simplification is pertinent in
that it can also be viewed as a particular case of a maybe more natural 
simplification of forest fire models.  
Roughly, fire forest model can be described as follows: on a lattice, trees
are born (sites become occupied) at a certain rate, say $1$ ; 
at each tree, a fire may start at 
some rate, say $\lambda>0$: the site becomes vacant and fire propagates to 
neighbouring trees (occupied sites) at a  given speed (see~\cite{drossel} 
for a precise description). Taking an infinite  propagation speed means that
the whole connected component (of sites occupied by trees) containing
the ignited tree burns at once (one may think of lightning). Our model
corresponds to the case $\lambda=1$, infinite propagation speed and a 
lattice equal to $\zz$. From
the point of view of SOC systems, the interesting phenomenon is in the
asymptotic regime  $\lambda
\to 0$. Indeed fires are less frequent, but when they occur, destroyed  
clusters may be huge.  These models have been subject to a lot of numerical 
and heuristical studies (see~\cite{grassberger} for references), but 
fewer rigorous results. Even existence and uniqueness 
of the process for a multidimensionnal lattice and given $\lambda$ 
has been proved only recently \cite{durre1,durre2}. 
Limiting rescaling when $\lambda \to 0$ has been studied numerically
\cite{drossel, grassberger} but attempts to give a rigorous
basis, even in dimension 1  are more recent \cite{vdbj,vdbb,brouwer}. 
Still our model had not received a complete rigorous
treatment, and as far as we understand, even if results are not
surprising they are now quite complete and the approach we propose 
may be extended.

\vip
  
Consider the  model in which {\it birth} flocks follow Poisson processes
with rate $1$, while {\it killing} flocks follow Poisson
processes with rate $\lambda>0$. We believe that our result could
be extended without difficulty to 
the case where $\lambda \geq 1$ (so that the clusters are
not very large). 
In the case where $\lambda<1$, the method we use probably 
breaks down, but the refined version of the algorithm described in 
Appendix \ref{appendix} gives 
hints for further research in this direction

\vip

The paper is organized as follows.
In Section \ref{coupling}, we show that the avalanche process
can be coupled with (and compared to) a very simple system of independent
particles which we call a Bernoulli process. The invariant distribution
of this particle system is an infinite product of Bernoulli distributions.

In Section \ref{reverse}, we show how to build the invariant distribution 
of the avalanche process from a stationnary Bernoulli process on an 
a.s.~finite time intervall, provided some cluster (concerning essentially
the Bernoulli process) is a.s.~finite.

We obtain some large-deviation type upperbounds for the width and height 
of this cluster 
in Section \ref{estimate}.

This allows us to conclude the proof in Section \ref{concl}:
the invariant distribution exists and can be perfectly
simulated. We can estimate 
the decay of correlations in the invariant distribution
of the avalanche process, using the upperbound of 
the width of the previously cited cluster. The coupling also shows, in some
sense, the uniqueness of the invariant distribution and
the trend to equilibrium. The rate of return to equilibrium
is obtained as a corollary of the upperbound of the height of the 
cluster.
In Appendix \ref{appendix}, we write down
the {\it perfect} simulation algorithm for the invariant distribution
derived from Sections \ref{reverse} and \ref{estimate}.

We finally introduce a related coagulation-fragmentation 
mean-field model in Section \ref{cf}: assuming that the correlations
between the sizes of connected components of occupied sites are neglictable,
we write down an infinite system of ordinary differential
equations satisfied by the {\it concentrations}
of clusters with size $k$, for $k\geq 1$: each pair of clusters
coalesce at constant rate, while each cluster with size $k$ breaks
into clusters with size $0$ at rate $k$.
The equilibrium state of the system of O.D.E.s can be computed 
almost explicitely. Numerical experiments
show that this model is an excellent approximation of
the avalanche process, at least from a {\it global} point of 
view.

\section{The coupling with a Bernoulli process}\label{coupling} 
\setcounter{equation}{0}

The starting point of our results is that we may deduce a 
realization of the (possible) equilibrium $\Pi$ of the avalanche process
from that of a much simpler process, which we now describe.

\vip

Consider as before an IFPP $N$, and an initial state 
$\zeta_0 \in \epr$.
Assume that the snow flocks are falling on each site $i$ according
to $N(i)$, but that the avalanche is restricted to the site $i$:
if $i$ was vacant, it becomes occupied as before, but if it was occupied, 
it becomes vacant, letting its neighbors enjoy their own life.
Denote, for each $i \in \zz$, each $t\geq 0$, by $\zeta_t(i)=1$
(resp. $\zeta_t(i)=0$)
if the site $i$ is occupied (resp. vacant) at time $t$.

The process $(\zeta_t)_{t \geq 0}$ is obviously well-defined, 
unique, and explicit: for $i\in\zz$, $t\geq 0$,
$\zeta_t(i)=\zeta_0(i)$ (resp. $\zeta_t(i)=1-\zeta_0(i)$)
if $N_t(i)$ is even (resp. odd). In other words,
\begin{equation}
\zeta_t(i)=\frac{1}{2}[1 - (-1)^{\zeta_0(i)+N_t(i)}].
\end{equation}
We call it the $\zeta_0$-Bernoulli process (or if necessary the
$(\zeta_0,N)$-Bernoulli process).


Let us
now describe its trend to equilibrium.
 
\begin{lem}\label{ttebp}
Let $\Gamma= \otimes_{i\in\zz}
\left(\frac{1}{2}\delta_0+\frac{1}{2}\delta_1 \right)$ be the infinite 
product of Bernoulli laws with parameter $1/2$.
For $\zeta_0 \in \epr$, denote by $\Gamma_t^{\zeta_0}$ 
the law of the $\zeta_0$-Bernoulli process at time $t$. Then for all
$l \geq 0$, all $t\geq 0$,
\begin{equation}
\sup_{\zeta_0\in\epr}
\left|(\Gamma_t^{\zeta_0})_{[-l,l]}- \Gamma_{[-l,l]} \right|_{TV} \leq 
(2l+1)e^{-2t}.
\end{equation}
As a consequence, $\Gamma$ is the only invariant distribution of
the Bernoulli process.
\end{lem}

\begin{proof}
Let thus $\zeta_0 \in \epr$ be fixed and $(\zeta_t)_{t\geq 0}$ be
the  $(\zeta_0,N)$-Bernoulli process, for some given IFPP $N$.
First of all observe that for any $t\geq 0$, any $i \in \zz$, 
using the explicit expression of $\zeta_t(i)$ leads us to
\begin{eqnarray}
P[\zeta_t(i)=0]&=&\indiq_{\{\zeta_0(i)=0\}}P[N_t(i) \hbox{ is even}]
+ \indiq_{\{\zeta_0(i)=1\}}P[N_t(i) \hbox{ is odd}] \ala
&=& \frac{1}{2}\left[\indiq_{\{\zeta_0(i)=0\}} (1+e^{-2t}) + 
\indiq_{\{\zeta_0(i)=1\}} (1-e^{-2t}) \right].
\end{eqnarray}
This implies that $|P[\zeta_t(i)=0]- \frac{1}{2}| \leq e^{-2t}/2$.
By the same way, $|P[\zeta_t(i)=1]- \frac{1}{2}| \leq e^{-2t}/2$, so that
we get $\left|(\Gamma_t^{\zeta_0})_{\{i\}}- 
\Gamma_{\{i\}} \right|_{TV} \leq e^{-2t}$. The result follows, since $[-l,l]$
contains $2l+1$ sites and since the coordinates of 
$(\zeta_t(i))_{i \in [-l,l]}$
are independent.
\end{proof}

Next, we explain how to reverse time in the stationnary 
Bernoulli process. This will be usefull to build {\it from the past} 
the invariant distribution of the avalanche process.

\begin{lem}\label{revtimeber}
Let $\zeta_0$ be a $\epr$-valued random variable with law $\Gamma$.
Consider $N$ an IFPP, and let $(\zeta_t)_{t\geq 0}$ be the (stationnary)
$(\zeta_0,N)$-Bernoulli process.

Consider the process $(\tzeta_t)_{t \in (-\infty,0]}$ defined
by $\tzeta_t=\zeta_{(-t)-}$ for all $t\leq 0$. Then this is again
a (stationnary) Bernoulli process, in the sense that 
for any $T \leq 0$, the process $(\tzeta_t)_{t\in [T,0]}$ 
is a $(\tzeta_T,N^T)$-Bernoulli process 
with $\tzeta_T$ independent of $N^T$, with $\tzeta_T\sim\Gamma$,
and where the IFPP $N^T$ on $[T,0]$ is defined
by $N^T_t(i)= N_{(-T)-}(i)-N_{(-t)-}(i)$ for $t\in [T,0]$ and $i \in \zz$.
\end{lem}

\begin{proof}
Let $T=-S<0$ be fixed.
Using the explicit formula, we know that for all $i\in\zz$, all $t\in[T,0]$,
\begin{eqnarray}
\tzeta_t(i)&=&\zeta_{(-t)-}(i)= 
\frac{1}{2}[1 - (-1)^{\zeta_0(i)+N_{(-t)-}(i)}]\ala
&=&\frac{1}{2}[1 - (-1)^{\zeta_0(i)+N_{(-T)-}(i)-N^T_t(i)}]
=\frac{1}{2}[1 - (-1)^{\tzeta_T(i)+N^T_t(i)}],
\end{eqnarray}
since one easily checks that $\zeta_0(i)+N_{(-T)-}(i)$ 
and $\tzeta_T(i)=\zeta_{(-T)-}(i)$
have the same parity.
Thus, we just have to prove that (i) $\tzeta_T\sim\Gamma$, (ii) $N^T$ 
is an IFPP on $[T,0]$, (iii)  $\tzeta_T$ and $N^T$ are independent.

Point (i) is obvious from the stationnarity of $(\zeta_t)_{t\geq 0}$,
while point (ii) is a well-known fact about Poisson processes. 
To prove point (iii), it suffices to notice that for any $i\in \zz$,
$x\in\{0,1\}$, 
\begin{eqnarray}
P\left[\tzeta_T(i)=x \left\vert \sigma((N^T_t(i))_{t\in[T,0]}\right. \right]
=P\left[\zeta_{S-}(i)=x \left\vert \sigma((N_t(i))_{t\in[0,S]} 
\right.\right]\ala
= \indiq_{\{N_{S-}(i) \hbox{ is even}\}}P[\zeta_0(i)=x] + 
\indiq_{\{N_{S-}(i) \hbox{ is odd}\}}P[\zeta_0(i)=1-x]\ala
=\frac{1}{2}= P\left[\tzeta_T(i)=x \right].
\end{eqnarray}
This ends the proof.
\end{proof}

We will also need later the following monotonicity  result about the Bernoulli
process.

\begin{lem}\label{bermono}
Consider $N$ and $V$ two independent IFPPs.
Let $\zeta_0^1, \zeta_0^2 \in \epr$.
Consider the $(\zeta_0^1,N)$-Bernoulli process
$(\zeta_t^1)_{t\geq 0}$.

There exists $M$ an IFPP such that, denoting by  
$(\zeta_t^2)_{t\geq 0}$ the $(\zeta_0^2,M)$-Bernoulli process, a.s.,
for all $t\geq 0$, all $i\in \zz$,

(i) $M_t(i) = \int_0^t \indiq_{\{\zeta_\sm^1(i)=\zeta_\sm^2(i)\}} dN_s(i) + 
\int_0^t \indiq_{\{\zeta_\sm^1(i) \ne \zeta_\sm^2(i)\}} d V_s(i)$;

(ii) if $\gamma_i:=\inf\{t \geq 0; \zeta^1_t(i)=\zeta^2_t(i)\}$,
$P[\gamma_i \geq t]\leq e^{-2t}$, and $(M_{\gamma_i+t}-M_{\gamma_i})_{t\geq 0}
= (N_{\gamma_i+t}-N_{\gamma_i})_{t\geq 0}$;

(iii) if $\zeta_t^1(i) = \zeta_t^2(i)$ then 
$\zeta_{t+s}^1(i) = \zeta_{t+s}^2(i)$ for all $s\geq 0$;

(iv) if $\zeta_0^1(i)\leq \zeta_0^2(i)$, then  
$\zeta_t^1(i) \leq \zeta_t^2(i)$ a.s. for all $t\geq 0$.

We will say that $(\zeta_t^1,\zeta^2_t)_{t\geq 0}$ are the 
$(\zeta_0^1,\zeta_0^2,N,V)$-coupled Bernoulli processes.
\end{lem}

Of course, the more natural coupling consisting in building the two
Bernoulli processes with the same IFPP would not preserve order
as time evolves.

\vip

\begin{proof}
The coupling we use here consists
in choosing the same Poisson process $N(i)$ for both processes when
$\zeta_0^1(i) = \zeta_0^2(i)$, so that they will appear or die 
simultaneously, and will remain equal for all times.
But if $0=\zeta_0^1(i) < \zeta_0^2(i)=1$ (resp. 
$0=\zeta_0^2(i) < \zeta_0^1(i)=1$), 
we use first independent Poisson processes:
$\zeta_t^2(i)$ dies 
using $V_t(i)$, while $\zeta_t^1(i)$ appears following
$N_t(i)$: after this first jump, they become equal, and we then use the same
Poisson process $N_t(i)$, and they remain equal for all times.

\vip

More rigorously,
for $i \in\zz$, denote by $T_i$ (resp. $S_i$) the first instant of
jump of $N(i)$ (resp. $V(i)$), and put $\tau_i=T_i \land S_i$.
It is immediate that $\tau_i$ follows an exponential distribution
with parameter $2$.
Define the process $M(i)$ by
\begin{equation}
M_t(i)= \indiq_{\{\zeta_0^1(i)=\zeta_0^2(i)\}} N_t(i)
+ \indiq_{\{\zeta_0^1(i)\ne \zeta_0^2(i)\}} 
\left[V_{t\land \tau_i}(i) + (N_t-N_{\tau_i})
\indiq_{\{t > \tau_i\}} \right].
\end{equation}
Then $M(i)$ is classically a Poisson process with rate $1$. We thus may
define the $(\zeta_0^2,M)$-Bernoulli process $(\zeta^2_t)_{t\geq 0}$.

Let us check points (i), (ii), (iii), and (iv). If $\zeta_0^1(i)=\zeta_0^2(i)$,
these points are obvious and $\gamma_i=0$, since then $M(i)=N(i)$ and 
$\zeta_t^1(i)=\zeta_t^2(i)$ for all times.

If $0=\zeta_0^1(i)<\zeta_0^2(i)=1$, then $\gamma_i=\tau_i$ and
$0=\zeta_t^1(i)<\zeta^2_t(i)=1$ for $t\in[0,\gamma_i)$. 
Easy considerations show that
for  $t\geq \gamma_i$,
$M_t(i)$ and $N_t(i)$ have an opposite parity, which implies
that $\zeta_t^1(i)=\zeta_t^2(i)$. This shows points (ii), (iii), and (iv).
Since $\gamma_i= \inf \{t\geq 0, \zeta^1_t(i)=\zeta^2_t(i)\}$,
point (i) can be written as $M_t(i)=V_{t\land \gamma_i} + (N_t-N_{\gamma_i})
\indiq_{\{t > \gamma_i\}}$, which achieves the proof.
\end{proof}

We now describe the coupling between the avalanche and Bernoulli processes.

\begin{prop}\label{avaber}
Consider $N$ and $V$ two independent IFPPs.
Let $\eta_0\in E$ and $\zeta_0 \in \epr$. 
Assume that for all $i \in \zz$, $\eta_0(i) \leq \zeta_0(i)$.
Consider the $(\zeta_0,N)$-Bernoulli process
$(\zeta_t)_{t\geq 0}$.

There exists $M$ an IFPP such that, denoting by  
$(\eta_t)_{t\geq 0}$ the $(\eta_0,M)$-avalanche process, a.s.,
for all $t\geq 0$, all $i\in \zz$,

(i) $\eta_t(i) \leq \zeta_t(i)$; 

(ii) $M_t(i) = \int_0^t \indiq_{\{\eta_\sm(i)=\zeta_\sm(i)\}} dN_s(i) + 
\int_0^t \indiq_{\{\eta_\sm(i)  < \zeta_\sm(i)\}} d V_s(i)$.

We will say that $(\zeta_t,\eta_t)_{t\geq 0}$ is the 
$(\zeta_0,\eta_0,N,V)$-coupled Bernoulli-avalanche process.
\end{prop}

Again here, building the Bernoulli and avalanche processes with the same
IFPP would not preserve the order.

\vip

\begin{proof}
The coupling is the following.
For each $i\in \zz$, at time $t\geq 0$, we use:

(a) the same IFPP $N_t(i)$ to make appear a flock in $\eta$ and $\zeta$  
if $\eta_{t-}(i)=\zeta_{t-}(i)=0$;

(b) the same IFPP $N_t(i)$ to make die the flock at $i$ 
(in $\zeta$) or the whole connected component of flocks around $i$ 
(in $\eta$) if $\eta_{t-}(i)=\zeta_{t-}(i)=1$;

(c) the IFPP $N_t(i)$ to make die the flock at $i$ (in $\zeta$)
and the independent IFPP $V_t(i)$ to make appear a flock (in $\eta$) if 
$0=\eta_{t-}(i)<\zeta_{t-}(i)=1$.

This construction guarantees that for all $t\geq 0$, all $i\in\zz$,
$\eta_t(i) \leq \zeta_t(i)$. The rigorous proof is similar to that
of Lemma \ref{bermono}.
\end{proof}

This coupling is illustrated by Figure 1, and can be represented
graphically in the following way.

\begin{gc}\label{gcab}
(a) Initially, each site of $\zz$ is occupied or not
according to $\zeta_0$ or $\eta_0$.
We draw black (resp. grey) segments to represent the marks of 
$N$ (resp. $V$) above each site of $\zz$.

(b) Next, we deduce the Bernoulli process $\zeta$: 

when an occupied site encounters a black mark,
it becomes vacant;

when a vacant site encounters a black mark,
it becomes occupied.

(The Bernoulli process is not concerned
with the grey marks).

(c) Finally, we deduce the avalanche process $\eta$: 

when an occupied site, say $i$, encounters a 
black mark, this makes become vacant the whole connected component 
of occupied sites around $i$;

when a vacant site (say the site $i$, at time $t$) 
encounters a black mark, it becomes 
occupied (and so does it in the process
$\zeta$) if and only if the Bernoulli process satisfies 
$\zeta_{t-}(i)=0$;

when a vacant site (say the site $i$, at time $t$), encounters a 
grey mark, then it becomes occupied if and only if $\zeta_{t-}(i)=1$.
\end{gc}

This graphical construction is possible because $\eta_0\in E$, 
which guarantees that for any $T>0$, there are a.s.~infinitely many
sites $i$ for which $\eta_0(i)=N_T(i)=V_T(i)=0$, and since
such sites {\it cut} the interactions.

\begin{figure}[t]
\centerline{Figure 1: Coupled avalanche and Bernoulli processes.}
\vskip0.1cm 
\centerline{\hskip1cm\includegraphics{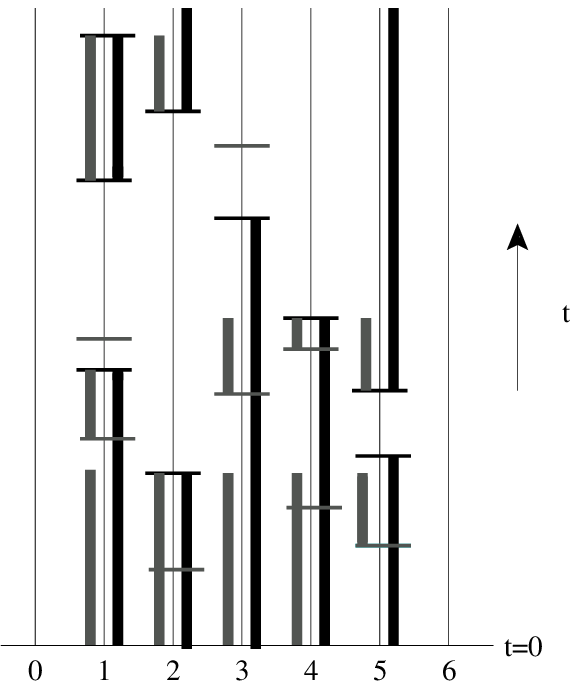}}
\begin{footnotesize}
The Bernoulli (resp. avalanche) process is represented in black on the right
(resp. in grey on the left) of each site.
Initially, the Bernoulli (resp. avalanche)
process is occupied on the sites $1,2,3,4,5$ (resp. $1,2,3,4$), 
and vacant on the sites $0,6$ (resp. $1,5,6$).
The Bernoulli process is easily constructed from the black marks:
each time a site encounters a black mark, its state changes.
Next, we have to build the avalanche process. 
The sites $2$ and $4$ are not affected
by the first grey marks, since they are occupied. 
On the contrary, $5$ becomes occupied when it encounters its first grey mark,
since it is vacant and the Bernoulli process is occupied 
(at this time on this site).
Next, the 
site $2$ encounters a black mark, which kills him and its whole connected
component of occupied sites, that is $1,2,3,4,5$. 
Next, the (vacant) site $5$ encounters a black mark, but it
does not become ocuupied, because the Bernoulli process is occupied. 
Next, the site
$3$ encounters a grey mark: since it is vacant and the Bernoulli
process is occupied, it becomes occupied. But it is killed 
again by the site $4$, which becomes vacant because it 
encounters a black mark. And so on...
\end{footnotesize}
\end{figure}

An immediate consequence of Proposition \ref{avaber} is the following,
which we will use later.

\begin{cor}\label{tropfastoche}
Let $\Pi$ be an invariant distribution of the avalanche process.
Recall that $\Gamma$ is the invariant distribution  of the Bernoulli
process.
Then $\Pi$ is stochastically smaller than  $\Gamma$.
This implies that for any random variable $\zeta\sim\Gamma$, 
we may find
a random variable  $\eta\sim\Pi$ such that a.s.,
for all $i\in\zz$, $\eta(i) \leq \zeta(i)$. 
\end{cor}

\begin{proof}
First, Supp $\Pi \subset E$, since the
rate of death for each site is bounded below by 1.
Consider $\eta_0=\zeta_0 \sim \Pi$.  Using Proposition \ref{avaber}, 
consider a $\eta_0$-avalanche process and a $\zeta_0$-Bernoulli 
process such that
a.s., for all $t\geq 0$, $i\in\zz$, $\eta_t(i) \leq \zeta_t(i)$. 
Of course, $\eta_t \sim \Pi$ for all $t$, while $\zeta_t$ goes in 
law to $\Gamma$
as $t$ tends to infinity, due to Lemma \ref{ttebp}.
We conclude that for any $\gamma \in \epr$, setting 
$F_\gamma=\{\alpha \in \epr,
\; \forall \; i\in \zz,\; \alpha(i)\geq \gamma(i)\}$,  
$\Pi(F_\gamma) \leq \Gamma(F_\gamma)$. This says exactly that 
$\Pi$ is stochastically smaller than $\Gamma$.
\end{proof}

\section{Coupling the invariant distributions}\label{reverse}
\setcounter{equation}{0}

Our aim in this section is to describe a way to build 
the invariant distribution of the avalanche process from
that of the Bernoulli process, using the coupling introduced in 
Proposition \ref{avaber}. Our method is based on the ideas of the famous
Propp-Wilson algorithm, \cite{pw}, which concerns Markov chains
with finite state space.
In the sequel, we will denote, for $\zeta \in \{0,1\}^\zz$,
\begin{equation}\label{dfez}
E_\zeta:= \left\{ \eta \in E; \forall \; 
i\in\zz, \; \eta(i) \leq \zeta(i) \right\}.
\end{equation}

\begin{prop}\label{ps1}
Let $V$ and $N$ be two independent IFPPs, and $\zeta_0\sim \Gamma$ 
(recall Lemma \ref{ttebp}).
Consider the $(\zeta_0,N)$-Bernoulli process $(\zeta_t)_{t\geq 0}$, 
and its time-reversed
$(\tzeta_t)_{t\in(-\infty,0]}$ built in Lemma \ref{revtimeber}.

For $T \in(-\infty,0]$ and $\varphi \in E_{\zeta_T}$, we denote by
$(\tzeta_t,\eta^{T,\varphi}_t)_{t\in [T,0]}$ the 
$(\tzeta_T,\varphi,N^T,V^T)$-coupled 
Bernoulli avalanche process 
with $N^T_t(i)=N_{(-T)-}(i)-N_{(-t)-}(i)$ and 
$V^T_t(i)=V_{(-T)-}(i)-V_{(-t)-}(i)$
for $t\in [T,0]$ and $i \in \zz$.

Observe that a.s., due to Proposition \ref{avaber}, 
for all $S\leq T \leq t\leq 0$, all $\varphi \in E_{\tzeta_S}$,
\begin{equation}\label{jaja}
\eta^{S,\varphi}_T \in E_{\tzeta_T} \hbox{ and } 
\eta^{S,\varphi}_t = \eta^{T,\eta^{S,\varphi}_T }_t.
\end{equation}

Denote, for each $i \in \zz$, by (here $\bZ\in E$ is the state with all sites
vacant)
\begin{equation}
\tau_i=\sup \{T \leq 0 ;
\forall \varphi\in E_{\tzeta_T},\; \eta^{T,\varphi}_0(i)=\eta^{T,\bZ}_0(i)\}.
\end{equation}
Assume for a moment that a.s., for all $i \in \zz$, $\tau_i>-\infty$. 
Notice that we have a.s., for all $i\in\zz$,
all $s_1\leq s_2<0$, all $\varphi_1 \in E_{\tzeta_{\tau_i+s_1}}$,
$\varphi_2\in E_{\tzeta_{\tau_i+s_2}}$,
\begin{equation}\label{jojo} 
\eta^{{\tau_i+s_1},\varphi_1}_0(i)=\eta^{{\tau_i+s_2},\varphi_2}_0(i).
\end{equation}
We thus may define $\eta_0(i)=\eta^{{\tau_i+s},\bZ}_0(i)$ (which
does not depend on $s<0$). 

Then $\Pi:=\cL(\eta_0)$ is the unique invariant distribution of
the avalanche process. 
\end{prop}

It seems that the Bernoulli process is almost unusefull in this statement.
However, it allows us to couple all the avalanche processes
(with different initial conditions) together.
Furthermore, the behaviour of $\tau_i$ will be studied through
the Bernoulli process. For example,
notice that $\tau_i=0>-\infty$ if $\zeta_0(i)=0$.
Indeed, due to Proposition \ref{avaber}, we know that for all $T<0$, 
all $\varphi \in E_{\tzeta_T}$, all $s\in [T,0]$, $\eta_t^{T,\varphi}(i) \leq 
\tzeta_t(i)$,
which implies that if $\tzeta_0(i)=0$ (i.e. $\zeta_0(i)=0$), then   
$ \eta_0^{T,\varphi}(i)=0$. Hence $\tau_i=0$ and $\eta_0(i)=0$.
When $\zeta_0(i)=1$, it is much less clear that $\tau_i>-\infty$.

\vip

\begin{proof} We split the proof into three parts.

{\bf Step 1.} Let us first explain (\ref{jaja}) and (\ref{jojo}). 
First, the fact that for $S\leq T\leq 0$ and $\varphi \in E_{\tzeta_S}$,
$\eta^{S,\varphi}_T \in E_{\tzeta_T}$ is straightforward from
Proposition \ref{avaber}. Then the second equality in (\ref{jaja})
follows from the construction.
Next, consider $i\in\zz$,
$s_1\leq s_2<0$, $\varphi_1 \in E_{\tzeta_{\tau_i+s_1}}$ and
$\varphi_2\in E_{\tzeta_{\tau_i+s_2}}$. Due to the definition
of $\tau_i$, we deduce that 
$\eta^{{\tau_i+s_1},\varphi_1}_0(i)=\eta^{{\tau_i+s_1},\bZ}_0(i)$.
On the other hand, we get from (\ref{jaja}) that 
$\eta^{{\tau_i+s_1},\bZ}_{\tau_i+s_2} \in E_{\tzeta_{\tau_i+s_2}}$ and
$\eta^{{\tau_i+s_1},\bZ}_0(i)
=\eta^{{\tau_i+s_2},\eta^{{\tau_i+s_1},\bZ}_{\tau_i+s_2}}_0(i)$,
the latter being equal to $\eta^{{\tau_i+s_2},\bZ}_0(i)$ due to
the definition of $\tau_i$. But using again the definition
of $\tau_i$, we deduce that $\eta^{{\tau_i+s_2},\varphi_2}_0(i)
=  \eta^{{\tau_i+s_2},\bZ}_0(i)$. This shows (\ref{jojo}).

\vip

{\bf Step 2.} Let us now show that $\Pi=\cL(\eta_0)$ is an 
invariant distribution for the avalanche process.
To this aim, call $(\eta^\bZ_t)_{t\geq 0}$ the $\bZ$-avalanche process.
Consider also
a bounded function $\Phi: E \mapsto \rr$ depending only on a finite number of
coordinates, say $\Phi(\eta)=\Phi((\eta(k))_{|k|\leq n})$ for some
$n\geq 0$. We will show
that $\lim_{T \to + \infty} E[\Phi(\eta^\bZ_T)]=E[\Phi(\eta_0)]$, which
classically suffices to conclude.

\vip

Consider now the processes coupled as in the statement.
First, $E[\Phi(\eta^\bZ_T)]=E[\Phi(\eta^{-T,\bZ}_{0})]$ for all $T\geq 0$.
Next, on the set $\Omega_n(T)=\{\forall \; |i|\leq n, \; \tau_i > -T \}$,
$\Phi(\eta^{-T,\bZ}_{0})=\Phi(\eta_{0})$ a.s.
Since $P[\Omega_n(T)]$ increases to $1$ as
$T$ increases to infinity (because a.s., 
$\tau_{-n}\lor ... \lor \tau_n >-\infty$),
we deduce that $\lim_{T\to+\infty} 
E[\Phi(\eta^{-T,\bZ}_{0})]=E[\Phi(\eta_0)]$. This concludes
the second step.

\vip

{\bf Step 3.} Consider another invariant distribution $\Pi'$ of the avalanche
process. Let $T \geq 0$, and consider, using Lemma 
\ref{tropfastoche}, a random
variable $\varphi_T\sim \Pi'$ such that $\varphi_T\in E_{\zeta_T}$ a.s.
Consider, as in Step 2, a bounded function $\Phi: E \mapsto \rr$ depending 
only on a finite number of coordinates, say 
$\Phi(\eta)=\Phi((\eta(k))_{|k|\leq n})$ for some
$n\geq 0$, and set $\Omega_n(T)=\{\forall \; |i|\leq n, \; \tau_i > -T \}$.
Then on $\Omega_n(T)$, $\Phi(\eta^{T,\varphi_T}_{0})=\Phi(\eta_{0})$.
On the other hand, $\eta^{T,\varphi_T}_{0}\sim \Pi'$, since $\Pi'$
is invariant. Using that $P[\Omega_n(T)]$ increases to $1$ 
as $T$ tends to $\infty$, 
we easily conclude
that $\int \Phi d \Pi' = E[\Phi(\eta_{0})]$. Thus $\Pi' = \Pi$.
\end{proof}

\section{The contour process}\label{estimate}
\setcounter{equation}{0}

Our aim in this section is to define and study a process
which will allow us to estimate $\tau_i$, for $i\in\zz$, and to 
bound the number of sites
involved in the construction of $\eta_0(i)$, in order to estimate
the decay of correlations.

The first idea is the following: consider the occupied
zone in $\zz \times [0,\infty)$ of the Bernoulli process. Clearly,
if this occupied zone has no infinite connected components, then
$\tau_i$ is finite a.s.~for all $i\in\zz$. 
Indeed, each site $i$ would then a.s.~be 
encompassed by a vacant zone of the Bernoulli process, which
implies that the avalanche process is also vacant,
and {\it cuts} the interaction in some sense, which would allow
us to build $\eta_0(i)$ from the stationnary process, using
the graphical construction \ref{gcab}.

But such a consideration would
probably lead to a {\it fat tail} estimate of the distribution of $\tau_i$,  
because we are in a critical case (the proportion of space 
occupied by the Bernoulli process is $1/2$). A way to overcome
this difficulty is to make use of the grey marks (recall Figure 1), which 
also give us some information about $\eta_0(i)$. 

Let us now define 
the left and right {\it contour} processes,
keeping in mind the coupling between stationnary measures built in Proposition
\ref{ps1}.

\begin{defin}\label{dfcp}
Let $\zeta_0  \in E$, and $N,V$ be two independent IFPPs.
We consider the $(\zeta_0,N)$-Bernoulli process 
$(\zeta_t)_{t\geq 0}$ and we introduce the filtration
$\cG_t=\sigma\{\zeta_0(i),N_s(i),V_s(i); \; s\in [0,t], i \in \zz\}$.    

For $i\in\zz$, we define the 
$(\zeta_0,N,V)$-right contour process $(R_t^i)_{t\geq 0}$ 
around $i$, with values in $\zz+\frac{1}{2} \cup \{\infty\}$
(see Figure 2 for an illustration)
by 
\begin{equation}
R_t^i=\sum_{n\geq 0} R^i_{T_n^i} \indiq_{\{t \in [T_n^i,T_{n+1}^i)\}},
\end{equation} 
where:

Initially, $T_0^i=0$, 
$R_0^i = \inf\{k \geq i, \; \zeta_0(k)=0\} - \frac{1}{2}$. 
For $n\geq 0$, 
\begin{equation}
T_{n+1}^i= \inf \{t > T_n^i, \Delta N_t(R_{T_n^i}^i+\frac{1}{2})+ 
\Delta N_t(R_{T_n^i}^i-\frac{1}{2}) +  
\Delta V_t(R_{T_n^i}^i-\frac{1}{2})>0  \}.
\end{equation}
Then

(a) if $\Delta N_{T_{n+1}^i}(R_{T_n^i}^i+\frac{1}{2})>0$, then 
$R_{T_{n+1}^i}^i= \inf\{k>R_{T_n^i}^i, \; \zeta_{T_{n+1}^i}(k)=0\}
-\frac{1}{2}$;
(b) if $\Delta N_{T_{n+1}^i}(R_{T_n^i}^i-\frac{1}{2})>0$, then 
$R_{T_{n+1}^i}^i= \sup\{k<R_{T_n^i}^i, \; \zeta_{T_{n+1}^i}(k)=1\}+
\frac{1}{2}$;
(c) if $\Delta V_{T_{n+1}^i}(R_{T_n^i}^i-\frac{1}{2})>0$, then

\hskip0.5cm (i) if 
$\zeta_{T_{n+1}^i}(R_{T_n^i}^i-\frac{3}{2})=1$, $R_{T_{n+1}^i}^i=R^i_{T_n^i}$,

\hskip0.5cm (ii)  if $\zeta_{T_{n+1}^i}(R_{T_n^i}^i-\frac{3}{2})=0$, then 
$R_{T_{n+1}^i}^i= \sup\{k<R_{T_n^i}^i, \; \zeta_{T_{n+1}^i}(k)=1\}
+\frac{1}{2}$.

The left contour process $(L_t^i)_{t\geq 0}$ around $i$ 
is defined symmetrically.
\end{defin}

Remark that the sequence $(T^i_n)_{n\geq 0}$ contains all the instants 
of jumps of $(R^i_t)_{t\geq 0}$, but it contains also fictitious jumps
(case (c)-(i)).
We explain how to build graphically these contour processes, as illustrated
by Figure 2.

\begin{gc}
Draw above each site $i\in \zz$ the marks of $N$ in black and
those of $V$ in grey. Draw in black the Bernoulli process corresponding
to a given initial data $\zeta_0$. 

A time $0$, the right contour process
$R^0_0$ lies on the left of the first vacant site of $\zeta_0$ on the
right of $0$ (e.g., if $\zeta_0(0)=\zeta_0(1)=1$ and $\zeta_0(2)=0$,
then $R^0_0=1.5$).

Next the dynamics of $R^0$ are the following: 

(a) each time it encounters a black mark on its right, it jumps
to the left of the first vacant site on its right;

(b) when it encounters
a black mark on its left, it jumps to the right of the 
first occupied site on its left;

(c) when it encounters a grey mark on 
its left (say that $R^0_{t-} = i+0.5$), and if $\zeta_{t-}(i-1)=0$, then  
it jumps to the right of the first occupied site on the left
of $i-1$.

The process $(L^0_t)_{t\geq 0}$ follows the same dynamics, 
permuting the roles of {\rm left} and {\rm right}.
\end{gc}

\begin{figure}[t]\label{figcontour}
\centerline{Figure 2: The contour processes $(R^0_t)_{t\geq 0}$ and 
$(L^0_t)_{t\geq 0}$ around $0$.}
\vskip0.1cm
\centerline{ \hskip1cm \includegraphics{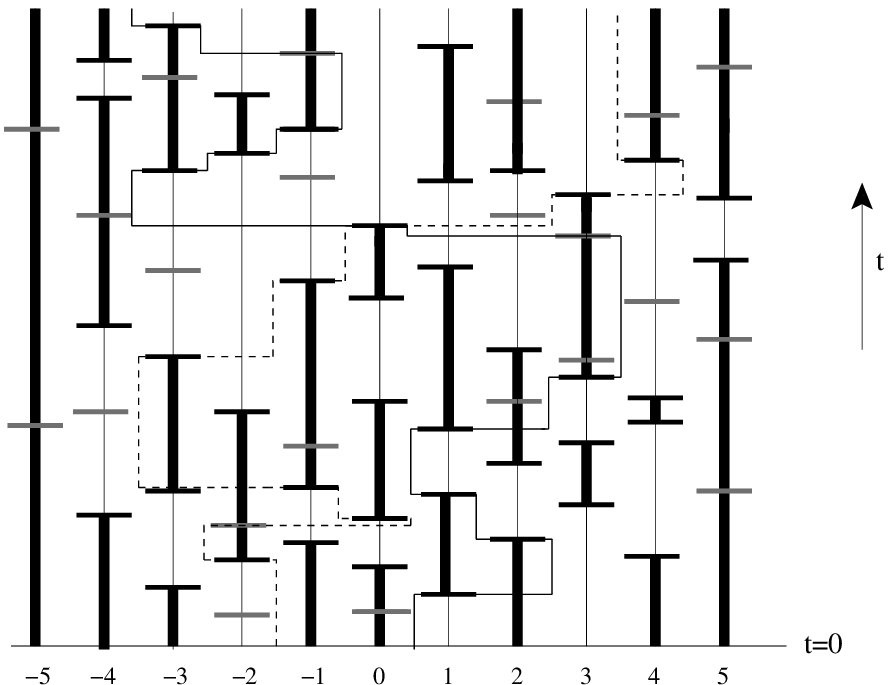} }
\begin{footnotesize}
The process $R^0_t$ (resp. $L^0_t$) is represented in plain (resp. dashed) 
line. First, $R^0_0=0.5$, since the first vacant site 
(of the Bernoulli process at time $0$) on the right of $0$
is $1$. By the same way, $L^0_0=-1.5$. Next, $R^0$ encounters a grey mark 
on its left, but since at this time $\zeta_t(-1)=1$, it does not jump.
Then $R^0$ encounters a black mark on its right, so that it jumps
to $2.5$, i.e. the left of $3$, which is (at this time) the first vacant
site on its right. Next, it encounters a black mark on its left and jumps 
to $1.5$, which is the right of $1$, i.e. the first vacant site on its left,
and so on... As we see on the picture, when it encouters its fourth
grey mark on its left, we have $R^0_{t-}=3.5$, and since $\zeta_{t-}(2)=0$,
it jumps to $0.5$, which is the right of $0$, i.e. the first occupied site
on the left of $2$.
\end{footnotesize}
\end{figure}

We will see in the next section that it is possible to build 
$\eta_0\sim \Pi$ (where $\Pi$ is the invariant distribution
of the avalanche process) in such a way that $\eta_0(i)$ depends only
on $\zeta_0$, $N$ and $V$ in the {\it box} delimited by $L^i_t$ and $R^i_t$
until they first meet (if they do). 
The main reason for this is the following property, which says
that in some sense, the contour processes encompass a given site
$i$ by a vacant zone of the Bernoulli process.

\begin{lem}\label{1firstpp}
We adopt the notations of Definition \ref{dfcp}.
A.s., for all $t\geq 0$, all $i\in \zz$, $\zeta_t(R^i_{t}+\frac{1}{2})= 
\zeta_t(L^i_{t}-\frac{1}{2})=0$.
\end{lem}

\begin{proof}
It is clear from the construction.
\end{proof}

Remark here that case (c)-(i) in Definition \ref{dfcp} is considered
to have this Lemma. Indeed, if we want the
right contour process to have only vacant sites (of the Bernoulli process)
on its right, we can use grey marks to jump to the left only when 
there is at least one vacant site on its strict left.

\vip

To study the decay of correlations,
we have to estimate the width of the {\it box}, while to study the rate
of trend to equilibrium, we have to estimate its height. 
The following estimates, central in our proof, 
will provide some bounds on these quantities.

\begin{prop}\label{est}
Let $\zeta_0 \sim \Gamma$, let $N,V$ be two independent IFPPs.
Consider the right and left $(\zeta_0,N,V)$-contour processes 
$(R_t^0)_{t\geq 0}$ 
and $(L_t^0)_{t\geq 0}$ around $0$. Consider the stopping time
(for the filtration $(\cG_t)_{t\geq 0}$)
\begin{equation}
\rho^0 = \inf\{t\geq 0; \; R_t^0 < L_t^0  \}
\end{equation}
and the random variable
\begin{equation}
\oR^0_\infty=\sup_{t\geq 0} R_t^0.
\end{equation}
(a) There exists $\beta>0$ such that $E[e^{\beta \rho^0}]<\infty$.

(b) There exists $\gamma>0$ such that $E[e^{\gamma \oR^0_\infty}]<\infty$.
\end{prop}

The remainder of this section is devoted to the proof of these
estimates. They seem quite natural, since
the process $(R^0_t)_{t\geq 0}$ is a sort of random walk with
negative {\it jump size expectation}: two types of events allow 
$(R^0_t)_{t\geq 0}$  to jump
to the left, while only one type allows him to jump to the right.
Furthermore, some symmetry seems to hold between the jumps to the 
right (due to $N$)
and those to left due to $N$.
Thus, the result seems almost obvious, and intuitively very clear.
However, we have not found a simple proof.
Of course, the main difficulty is that
$(R^0_t)_{t\geq 0}$  is not a continuous-time random
walk: its sizes of jumps are not independent. Thus quite a precise study
has to be done. Our strategy consists in bounding from above 
$(R^0_t)_{t\geq 0}$ by a continuous-time random walk with negative
jump size expectation.
We first describe some immediate 
properties of the contour processes.

\begin{lem}\label{firstpp}
We adopt the notations of Definition \ref{dfcp}.
Let $i \in \zz$ be fixed.

(a) If $\zeta_0\sim \Gamma$, then
the processes $(i-L^i_t)_{t\geq 0}$, $(R^i_t-i)_{t\geq 0}$, and
$(R^0_t)_{t\geq 0}$ have the same law (but are far from being independent).

(b) If $R_0^i<\infty$, then $R_t^i<\infty$ for all $t\geq 0$ a.s.

(c) For $j\leq i$, $R_t^j \leq R_t^i$ for all $t\geq 0$ a.s. 
Furthermore, if $R_t^j = R_t^i$ for some $t$, then a.s.,  
$R_{t+s}^j = R_{t+s}^i$ for all $s\geq 0$. 

(d) The counting processes
\begin{eqnarray}
Z^{1,i}_t= \sum_{n\geq 1} \indiq_{\{t \geq T_n^{1,i}\}} := 
\sum_{n\geq 1} \indiq_{\{t \geq T_n^i\}} \indiq_{\{\Delta
N_{T_{n}^i}(R_{T_{n-1}^i}^i+\frac{1}{2})>0\}}=N_t(R_{t-}^i+\frac{1}{2})  \ala
Z^{2,i}_t=\sum_{n\geq 1} \indiq_{\{t \geq T_n^{2,i}\}} := 
\sum_{n\geq 1} \indiq_{\{t \geq T_n^i\}} \indiq_{\{\Delta
N_{T_{n}^i}(R_{T_{n-1}^i}^i-\frac{1}{2})>0\}}=N_t(R_{t-}^i-\frac{1}{2})\ala
Z^{3,i}_t=\sum_{n\geq 1} \indiq_{\{t \geq T_n^{3,i}\}} := 
\sum_{n\geq 1} \indiq_{\{t \geq T_n^i\}} \indiq_{\{\Delta
V_{T_{n}^i}(R_{T_{n-1}^i}^i-\frac{1}{2})>0\}}=V_t(R_{t-}^i-\frac{1}{2})
\end{eqnarray}
are three independent Poisson processes with rate $1$. They are
$(\cG_t)_{t\geq 0}$-adapted, and independent of $\zeta_0$.
\end{lem}

Remark here that $Z^{1,i}$ counts the jumps to the right of $R^i$,
while $Z^{2,i}$ counts its jumps to the left due to $N$ 
(black marks on Figure 2)
and $Z^{3,i}$ counts its {\it possible} 
jumps to the left due to $V$ (grey marks on Figure 2).

\vip
 
\begin{proof}
Point (a) is obvious by symmetry and invariance by translation.
Point (b) follows from the fact that the Bernoulli a.s.~process belongs to 
$E$ for all $t>0$, even if it does not at time $0$.
Point (c) is clear from the construction.
Point (e) follows from classical properties on Poisson processes.
\end{proof}

We carry on with a natural monotonicity property.

\begin{lem}\label{toctoc}
We consider three independent IFPPs $N$, $V$ and $W$.
Let also $\zeta_0^1,\zeta_0^2 \in \epr$ satisfy, for all $i\in\zz$,
$\zeta_0^1(i) \leq \zeta_0^2(i)$. Then we build, recalling
Lemma \ref{bermono}, the $(\zeta_0^1,\zeta_0^2,N,W)$-coupled Bernoulli 
process $(\zeta^1_t,\zeta^2_t)_{t\geq 0}$. As stated in Lemma  
\ref{bermono}, $(\zeta^1_t)_{t\geq 0}$ is the 
$(\zeta_0^1,N)$-Bernoulli 
process, while  $(\zeta^2_t)_{t\geq 0}$ is the  
$(\zeta_0^2,M)$-Bernoulli process for some IFPP $M$.
We denote by $(R^{0,1}_t)_{t\geq 0}$ (resp. $(R^{0,2}_t)_{t\geq 0}$)
the $(\zeta_0^1,N,V)$ (resp. $(\zeta_0^2,M,V)$) right contour process around 
$0$.

We will say that $(R^{0,1}_t,R^{0,2}_t)_{t\geq 0}$ are the 
$(\zeta_0^1,\zeta_0^2,N,W,V)$-coupled
right contour processes around $0$. 
We have a.s., for all $t\geq 0$, $R^{0,1}_t \leq R^{0,2}_t$.
\end{lem}

\begin{proof}
The proof is obvious from the definition of the contour process,
since we know from Lemma \ref{bermono} that a.s., for all $t\geq 0$, all 
$i\in \zz$, $\zeta^1_t(i) \leq \zeta^2_t(i)$.
\end{proof}

We consider the following initial condition.

\begin{nota}\label{dfxi}
We say that a $\epr$-valued random variable $\tzeta_0$ 
has the distibution $\Xi$ is 
if $\tzeta_0(i)=1$ for $i\leq 0$, $\tzeta_0(1)=0$, and if
$(\tzeta_0(i))_{i\geq 2}$ is a family of i.i.d. Bernoulli random variables
with parameter $1/2$.
\end{nota}

Let us now explain our strategy to bound the right contour process
by a random walk:

we will first upperbound the initial configuration 
$\zeta_0$ of the Bernoulli process by a (possibly shifted) 
realization $\tzeta_0$ of $\Xi$,
we thus upperbound our contour process by the corresponding contour
process $\tR^0$,

then we will wait for the first instant $\tT_1^{1,0}$ of jump to the right 
of $\tR_t^0$ to the right; this yields a total jump which we will call
$Y_1:=\tR_{\tT_1^{1,0}}-\tR_0$, 
and whose expectation will be shown to be negative;

we will also observe that at $\tT_1^{1,0}$, 
we may again bound the configuration
of the Bernoulli process by a realization $\tzeta_0^1$ of $\Xi$ (shifted
around $\tR_{\tT_1^{1,0}}$) independent of $Y_1$;

this last renewal argument allows us to build, recursively,
a random walk with negative mean jump size, 
bounding from above our contour process.

\begin{lem}\label{renew}
Let $\tzeta_0 \sim \Xi$, and consider two independent IFPPs $N,V$. 
Consider the $(\tzeta_0,N,V)$-right contour process 
$(\tilde R_t^0)_{t\geq 0}$ around $0$, observe that $\tR^0_0=1/2$, and
denote by 
$\tT^{1,0}_1:= \inf\{t>0, \Delta \tR^0_t>0 \}$ the first instant 
where it jumps to the right. We also consider
$(\tzeta_t)_{t\geq 0}$ the $(\tzeta_0,N)$-Bernoulli process. 
We set $Y_1=\tR^0_{\tT^{1,0}_1}-1/2$. 

(i) Then $E[Y_1]<0$.

(ii) For all $\e\in (0,\ln 2)$, $E[e^{\e Y_1}]<\infty$.

Furthermore, then there exists $\tzeta_0^1\sim \Xi$ such that

(iii) a.s., $\tzeta_{\tT^{1,0}_1}(\tR^0_{\tT^{1,0}_1}+i-1/2) \leq 
\tzeta_0^1(i)$ for all $i \in \zz$,

(iv) $\tzeta_0^1$ is independent of $\cH_{\tT^{1,0}_1}$, where
$\cH_t = \sigma \{\tR^0_s, s \leq t \}$.
\end{lem}

\begin{proof}
To simplify the notation, we omit the superscript $0$ 
(which says that we are dealing with the contour process around $0$) 
in this proof.
We consider the three independent
Poisson processes with rate $1$ (see Lemma \ref{firstpp}-(d))
$\tZ^1_t= N_t(\tR_{t-} +1/2)$, $\tZ^2_t= N_t(\tR_{t-} - 1/2)$
and $\tZ^3_t= V_t(\tR_{t-} -1/2)$, and we denote by $(\tT^1_i)_{i\geq 1}$,
$(\tT^2_i)_{i\geq 1}$, $(\tT^3_i)_{i\geq 1}$, respectively, their 
successive instants of jumps. We also 
denote by $A_j=\cup_{i\geq 1} \{\tT^j_i\}$, for
$j=1,2,3$. We set $\tZ_t=\tZ^1_t+\tZ^2_t+\tZ^3_t$,
which is a Poisson process with rate $3$, we denote by
$(\tT_i)_{i\geq 1}$ its successive instants of jumps, and we set
$A=\cup_{i\geq 1} \{\tT_i\}=A_1 \cup A_2\cup A_3$. 
Finally, we also set for convenience $\tT_0=\tT^1_0=\tT^2_0=\tT^3_0=0$.
Recall that we want to study $Y_1=\tR_{\tT^1_1} - 1/2$.

\vip

{\bf Step 1.} For $n \geq 1$, the event
$\Omega_n=\{\tT_1 \notin A_1,...,\tT_{n-1} \notin A_1,\tT_n \in A_n\}$ 
occurs with probability $p_n:=\frac{2^{n-1}}{3^n}$, since $A_1,A_2,A_3$ are
the sets of jumps of three independent Poisson processes
with same rate. Notice also that on $\Omega_n$, $\tT^1_1=\tT_n$,
and we may write $Y_1= - \sum_{i=1}^{n-1} X_i + X_n$, where
$X_1,...,X_n$ are the successive sizes of the (possibly fictitious)
jumps of $\tR$ (at the instants $\tT_1<...<\tT_n$), with
$X_1\geq 0$, ..., $X_n \geq 0$. We obtain
\begin{equation}\label{xx1}
E[Y_1]= \sum_{n \geq 1} E[ \{-(X_1+...+X_{n-1}) + X_n\} \indiq_{\Omega_n}].
\end{equation}

\vip

{\bf Step 2.} 
Let us now bound from below  $C_{i,n}:= E[ X_i\indiq_{\Omega_n}]$, for
$1\leq i \leq n-1$. 

We denote by $Z_i$ the number of vacant sites of the Bernoulli process
on the strict left of $R_{\tT_i-}-1/2$ at time $T_i-$, 
that is 
\begin{equation}
Z_i:= \tR_{\tT_{i}-} -3/2 - \sup\{j \leq \tR_{\tT_{i}-}-3/2,\;
\tzeta_{\tT_i-}(j)=1\}.
\end{equation}
Then, due to the definition of $\tR$, we know that

(a) on $F^2_{i}:=\{\tT_i \in A_2\}$, 
$X_i=- \Delta \tR_{\tT_i} = 1+Z_i$,

(b) on $F^3_i:=\{\tT_i \in A_3\}$, 
$X_i=- \Delta \tR_{\tT_i} =  (1+Z_i)\indiq_{\{Z_i\geq 1\}}$.

Observe that $P[F^2_{i} \vert \Omega_n]=P[F^3_i \vert \Omega_n]=1/2$, and
that $F_i^2,F^3_i$ are independent of
$(Z_i,\tT_i)$ conditionally to $\Omega_n$. 
These are standard properties of Poisson processes. Hence,
\begin{eqnarray}\label{jjaabb}
C_{i,n}&=&\frac{1}{2}E\left[(1+Z_i)\indiq_{\Omega_n}+
(1+Z_i)\indiq_{\{Z_i\geq 1\}}\indiq_{\Omega_n} \right] \ala
&=& E\left[(1+Z_i)\indiq_{\Omega_n}\right] - \frac{1}{2}P[Z_i=0, \Omega_n]\ala
&=& P[\Omega_n]+ \sum_{k \geq 1}  P[Z_i \geq k, \Omega_n] 
- \frac{1}{2}P[Z_i=0, \Omega_n].
\end{eqnarray}
Let now $k \geq 1$ be fixed. We have
\begin{eqnarray}
P[Z_i \geq k, \Omega_n]=P\Big(
\tzeta_{\tT_i-}( \tR_{\tT_{i}-}-\frac{3}{2})=0,...,
\tzeta_{\tT_i-}(\tR_{\tT_{i}-}-\frac{1}{2}-k)=0, \Omega_n\Big) \ala
= E\left[ \prod_{l=1}^k P \left(\left.
\tzeta_{\tT_i-}(\tR_{\tT_{i}-}-\frac{1}{2}-l)=0  
\right\vert \Omega_n, \tT_i,\tT_{i-1} \right) \indiq_{\Omega_n} \right].
\end{eqnarray}
Indeed, recalling that on $\Omega_n$, $\tR$ has had only jumps 
to the left before $\tT_i$, 
we easily deduce that on $\Omega_n$,
the values of the Bernoulli process at sites 
$j \leq \tR_{\tT_{i}-}-\frac{3}{2}$ are mutually independent conditionnally
to $\tT_i,\tT_{i-1}$. 

Let us set $p_s:=(1-e^{-2s})/2=P[N_s \in 2\nn +1]$ for $s\geq 0$
(for $(N_t)_{t\geq 0}$ is a standard Poisson process with rate $1$).

Now for $l \geq 2$, the site $\tR_{\tT_{i}-}-\frac{1}{2}-l$
was occupied at time $0$, and its evolution is obviously independent
of $(\tR_t)_{t \in [0,\tT_i)}$, so that 
\begin{equation}
P \left(\left. \tzeta_{\tT_i-}(\tR_{\tT_{i}-}-\frac{1}{2}-l)=0  
\right\vert \Omega_n, \tT_i,\tT_{i-1} \right)= p_{\tT_i}.
\end{equation}

Next, the same argument holds for $l=1$ on the set $\{X_{i-1}>0\}$,
which indicates that the previous jump to the left
was not fictitious: we have 
\begin{equation}
P \left(\left. \tzeta_{\tT_i-}(\tR_{\tT_{i}-}-\frac{3}{2})=0  
\right\vert \Omega_n, \tT_i, \tT_{i-1},X_{i-1}>0 \right)= p_{\tT_i}.
\end{equation}
But on the event $\{X_{i-1}=0\}$, we know that 
$\tzeta_{\tT_{i-1}}(\tR_{\tT_{i}-}-\frac{3}{2})=1$. Hence we get
\begin{equation}
P \left(\left. \tzeta_{\tT_i-}(\tR_{\tT_{i}-}-\frac{3}{2})=0  
\right\vert \Omega_n, \tT_i, \tT_{i-1}, X_{i-1}=0 \right)= 
p_{\tT_i-\tT_{i-1}}.
\end{equation}
Noting that $p_{\tT_i} \geq p_{\tT_i-\tT_{i-1}}$, that $\{X_{i-1}=0\}
\subset \{\tT_{i-1} \in A_3\}$, we deduce that
\begin{equation}
P \left(\left. \tzeta_{\tT_i-}(\tR_{\tT_{i}-}-\frac{3}{2})=0  
\right\vert \Omega_n, \tT_i, \tT_{i-1}\right)
\geq \indiq_{\{\tT_{i-1} \in A_2\}}p_{\tT_i}+  
\indiq_{\{\tT_{i-1} \in A_3\}}p_{\tT_i-\tT_{i-1}}.
\end{equation}
Gathering the estimates obtained for $l\geq 2$ and $l=1$, we obtain,
for $k \geq 1$,
\begin{equation}
P[Z_i\geq k,\Omega_n]\geq E\left[\indiq_{\Omega_n} 
p_{\tT_i}^{k-1} \left( \indiq_{\{\tT_{i-1} \in A_2\}}p_{\tT_i}+  
\indiq_{\{\tT_{i-1} \in A_3\}}p_{\tT_i-\tT_{i-1}}  \right)\right]
\end{equation}

Using finally classical properties of Poisson processes, we see that
$\tT_{i-1},\tT_i$ are independent of $\Omega_n,\{\tT_{i-1}\in A_2\},
\{\tT_{i-1}\in A_3\}$ and that $P[\Omega_n \cap \{\tT_{i-1}\in A_2\}]
=P[\Omega_n \cap \{\tT_{i-1}\in A_3\}]=P[\Omega_n]/2=2^{n-1}/2.3^n$, so that
\begin{equation}
P[Z_i \geq k, \Omega_n] \geq \frac{2^{n-1}}{2.3^n} E\left[
\left(p_{\tT_i}+ p_{\tT_i - \tT_{i-1}} \right)
p_{\tT_i}^{k-1} \right].
\end{equation}
Next,
\begin{equation} 
P[Z_i=0 , \Omega_n] =  P[\Omega_n] -P[Z_i\geq 1 , \Omega_n]
\leq  \frac{2^{n-1}}{2.3^n} E \left[2-p_{\tT_i}-p_{\tT_i - \tT_{i-1}} \right].
\end{equation}
Thus, recalling (\ref{jjaabb}), for any $1\leq i \leq n-1$,
\begin{eqnarray}\label{xx2}
C_{i,n}
&\geq &  \frac{2^{n-1}}{3^n} E \left[1+\frac{p_{\tT_i}+ p_{\tT_i - \tT_{i-1}}}
{2-2p_{\tT_i}} -\frac{2-p_{\tT_i}-p_{\tT_i - \tT_{i-1}} }{4} \right] \ala
&\geq & \frac{2^{n-1}}{2.3^n}E\left[1+ \frac{p_{\tT_i}}{1-p_{\tT_i}}
+ \frac{p_{\tT_i - \tT_{i-1}}}{1-p_{\tT_i}} + \frac{1}{2}p_{\tT_i}+
\frac{1}{2}p_{\tT_i-\tT_{i-1}}\right]\ala
&=:& \frac{2^{n-1}}{2.3^n} (1+B_i),
\end{eqnarray}
where the last equality stands for a definition.

\vip

{\bf Step 3.}
We now upperbound $C_{n,n}:=E[X_n\indiq_{\Omega_n}]$. 
We denote by $Z_n$ the number of occupied sites on the strict right
of $\tR_{\tT_n-}$, that is
\begin{equation}
Z_n=\inf \{j \geq \tR_{\tT_n-}+3/2, \tzeta_{\tT_n-}(j)=0 \}-\tR_{\tT_n-}-3/2.
\end{equation}
By construction, we have $X_n= 1+Z_n$ on $\Omega_n$.
For $k \geq 1$, we set $J_k:=\tR_{\tT_n-}+k+ \frac{1}{2}$ and
$\xi_k:= \tzeta_{\tT_n-}(J_k)$. By construction, we have $X_n= 1+Z_n$ on 
$\Omega_n$, so that 
for $k \geq 1$
\begin{equation}
P[Z_n \geq k, \Omega_n] = P\left[\xi_1=1,...,\xi_k=1, \Omega_n \right].
\end{equation}
We now introduce the $\sigma$-field generated by the path of 
$(\tR_t)_{t \in [0,\tT_1^1)}$, containing also the fictitious jumps, that is,
for $\nu$ defined by $\tT_\nu=\tT_1^1$ ($\nu=n$ on $\Omega_n$),
\begin{equation}
\cH:= \sigma \left(\nu,\tT_1,...,\tT_\nu,
\Delta \tR_{\tT_1},...,\Delta \tR_{\tT_{\nu-1}} \right).
\end{equation}
Observe that obvisouly, $(J_k)_{k\geq 1}$ and $\Omega_n$ are $\cH$-measurable.

We will show that conditionnally to $\cH$, the sequence $(\xi_k)_{k\geq 1}$
is a family of independent random variables on $\Omega_n$, 
and that for each $k\geq 1$, 
$P[\xi_k=1 \vert \cH,\Omega_n]\leq 1/2$. 
Since $\Omega_n$ belongs to $\cH$,  for all $n\geq 1$, we will
deduce that 
\begin{equation}\label{jablat}
P[Z_n \geq k,\Omega_n]\leq \frac{1}{2^{k}}P[\Omega_n]=\frac{1}{2^{k}}
\frac{2^{n-1}}{3^n},
\end{equation}
so that (since $X_n=1+Z_n$ on $\Omega_n$),
\begin{equation}\label{xx3}
C_{n,n}=\sum_{k\geq 1} P[1+Z_n \geq k,\Omega_n]
=\sum_{k\geq 0} P[Z_n \geq k,\Omega_n]
\leq \frac{2^{n}}{3^n} .
\end{equation}

\begin{figure}[t]\label{figproof}
\centerline{Figure 3: Illustration of Step 3 (and Step 6).}
\vskip0.1cm
\centerline{ \hskip1cm \includegraphics{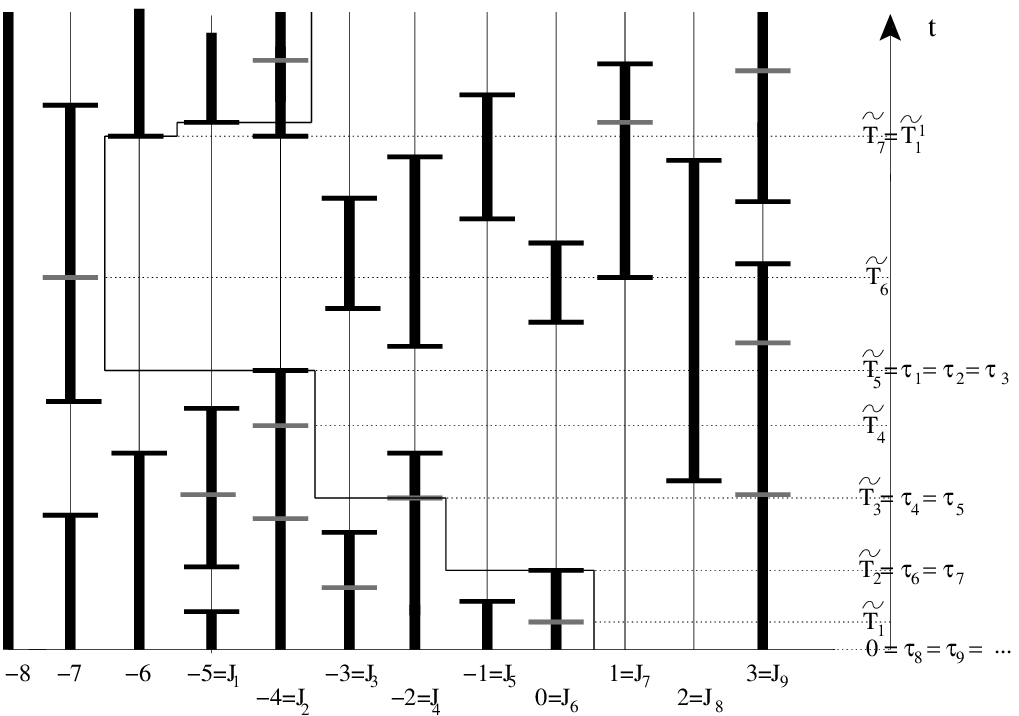} }
\begin{footnotesize}
With this realization, we have $\tR_{\tT_1^1-}=-6.5$ and 
$G=\{2,4,6\}$. We remark that the only site which is occupied 
by the Bernoulli process when
crossed by $\tR$ is the site $J_4=-2$: it is crossed through
a grey mark.
\end{footnotesize}
\end{figure}

Let us thus check the announced properties of the sequence $(\xi_k)_{k\geq 1}$.
For each $k\geq 1$, let $\tau_k=0$ if $J_k\geq 2$, and 
let $\tau_k$ be the unique instant $t\in [0,\tT^1_1)$ such that
$J_k-1 \in \{\tR_t,\tR_{t-}\}$ if $J_k\leq 1$. We refer to Figure 3 
for an illustration.
Roughly, $\tau_k$ is the last instant before $\tT_1^1$ where we get 
some information (from $\cH$) about the site $J_k$.
Of course, $(\tau_k)_{k\geq 1}$ is $\cH$-measurable.

\vip

We will show that the family of random variables $\tzeta_{\tau_k}(J_k)$
is mutually independent on $\Omega_n$ 
conditionnally to $\cH$, and that for each $k \geq 1$,
$P[\tzeta_{\tau_k}(J_k)=1\vert\cH, \Omega_n]\leq 1/2$, This will imply 
the announced properties, for two reasons:

(i) conditionnally to $\cH$ and $\Omega_n$, the evolution of the 
Bernoulli process at two different sites
$J_k$ (on $[\tau_k,\tT^1_1)$) and $J_l$ (on $[\tau_l,\tT^1_1)$)
are independent, 
since they concern independent Poisson processes,

(ii) for each $k\geq 1$, $P[\tzeta_{\tT_1^1-}(J_k)=1 \vert \cH,\Omega_n]
=P[\tzeta_{\tau_k}(J_k)=1\vert \cH,\Omega_n] (1-p_{\tT_1^1-\tau_k}) +
P[\tzeta_{\tau_k}(J_k)=0\vert \cH,\Omega_n] p_{\tT_1^1-\tau_k} \leq 1/2$.
Indeed, recall that $p_s=(1-e^{-2s})/2\leq 1/2$ stands for the probability that
a standard Poisson process at time $s$ is odd, and that
for $a,b \in [0,1/2]$, $a(1-b)+(1-a)b \leq 1/2$.

\vip

Consider now the random set, measurable with respect to $\cH$,
\begin{equation}
G:= \{ k\geq 1; \tau_{k+1}>\tau_{k+2}\}=\{ k\geq 1; J_k=\tR_{\tau_k-}-1/2\}.
\end{equation}
We notice that $\{J_k,k\in G\} \subset [\tR_{\tT_1^1},0]$.

We observe (see Figure 3) that conditionnally to $\cH$, $\Omega_n$,
for all $k \geq 1$, we have,
\begin{equation}\label{supertop}
\tzeta_{\tau_k}(J_k)= \indiq_{\{J_k \leq 1\}} \indiq_{\{k \in G\}} 
\indiq_{\{\tau_k \in A_3\}} + \indiq_{\{J_k \geq 2\}} \tzeta_0(J_k).
\end{equation}
Indeed, if $J_k\geq 2$, the formula is obvious because then $\tau_k=0$.
If $J_k \leq 1$, this comes from the fact that the only way for $\tR$
to jump to the left through an occupied site is that the jump
follows from a grey mark (i.e. $\tau_k \in A_3$) and that the concerned
site is just on the left of $\tR_{\tau_k -}$ (i.e. $k \in G$).

\vip

Conditionnally to $\cH$, $\Omega_n$,
the sequence of events $\{\tau_k \in A_3\}_{k \in G}$ is independent,
this assertion makes sense since $G$ is itself $\cH$-measurable.
Indeed, we know that conditionnally to $\cH$,
for $k \in G$,  $\{\tau_k \in A_3\}$ depends only on the Poisson
processes $N_t(J_k),V_t(J_k)$ (and on $\tzeta_{\tau_k-}(J_k-3/2)$ which
is $\cH$-measurable). The conditionnal independence (with respect to $\cH$,
$\Omega_n$)
of the family  $\{\tau_k \in A_3\}_{k \in G}$ follows then from
the fact that for $k_1<k_2$ in $G$, $J_{k_1}>J_{k_2}$ (see Figure 3), 
and from the 
independence of the Poisson processes $(N_t(i),V_t(i))$ and
$(N_t(j),V_t(j))$ for $i \ne j$.

\vip

Recall now (\ref{supertop}). 
Using the conditionnal independence (with respect to $\cH$, $\Omega_n$)
of the family  $\{\tau_k \in A_3\}_{k \in G}$, the 
fact that the sequence $(J_k)_{k\geq 1}$ is $\cH$-measurable,
and the fact that the family 
$(\tzeta_0(i))_{i \geq 2}$ is mutually independent and independent of
$\cH$, $\Omega_n$ (these are i.i.d. Bernoulli random variables 
with parameter $1/2$), we obtain that 
conditionnally to $\cH$, $\Omega_n$, the sequence 
$(\tzeta_{\tau_k}(J_k))_{k\geq 1}$ is mutually independent.

\vip 

We finally conclude by noting that for any $k\geq 1$,

\begin{eqnarray}
P\left[\left. \tzeta_{\tau_k}(J_k)=1 \right\vert \cH, \Omega_n\right]=
\indiq_{\{J_k\leq 1\}}\indiq_{\{k\in G\}} 
P \left[\left.\tau_k \in A_3  \right\vert \cH, \Omega_n \right]\ala
+ \indiq_{\{J_k\geq 2\}} P \left[\left.\tzeta_0(J_k)=1  
\right\vert \cH, \Omega_n \right] \leq 1/2.
\end{eqnarray}
The last inequality comes from the fact that 
$P \left[\tzeta_0(J_k)=1 \vert \cH, \Omega_n \right]=1/2$ if $J_k\geq 2$ 
as was previously noticed, while for $k \in G$, 
$P \left[\tau_k \in A_3 \vert \cH, \Omega_n \right] \leq 1/2$. Indeed,
having a look at Figure 3, we realize that for $k \in G$, 
$\{\tau_k \in A_3\} \subset \{\tzeta_{\tau_k-}(J_k-3/2)=0\}$, and that
due to classical properties of Poisson processes,
 \begin{eqnarray}
P \left[\tau_k \in A_3 \vert \cH, \Omega_n, 
\{\tzeta_{\tau_k-}(J_k-3/2)=0\}\right]= 1/2.
\end{eqnarray}

\vip

{\bf Step 4.} 
Gathering (\ref{xx1}), (\ref{xx2}) and (\ref{xx3}), we get
\begin{eqnarray}
E[Y_1]&=& \sum_{n \geq 1} \left\{C_{n,n} - 
\sum_{i=1}^{n-1} C_{i,n} \right\}\ala
&\leq&\sum_{n \geq 1} \left\{ \frac{2^{n}}{3^n}  - \sum_{i=1}^{n-1}
\frac{2^{n-1}}{2.3^n} (1+B_i) \right\} 
=2 - \sum_{i\geq 1} \frac{2^{i-1}}{3^i} (1+B_i)\ala
&\leq & 1 - \sum_{i\geq 1} \frac{2^{i-1}}{3^i} B_i.
\end{eqnarray}
To conclude that $E[Y_1]<0$, we thus have to prove that 
$I=\sum_{i\geq 1} \frac{2^{i-1}}{3^i} B_i >1$. But, we may
write $I=I_1+I_2+I_3/2+I_4/2$, with
\begin{eqnarray}
I_1:=\sum_{i\geq 1} \frac{2^{i-1}}{3^i} 
E \left[\frac{p_{\tT_i}}{1-p_{\tT_i}} \right], &&
I_2:=\sum_{i\geq 1} \frac{2^{i-1}}{3^i} 
E \left[\frac{p_{\tT_i - \tT_{i-1}}}{1-p_{\tT_i}} \right] \ala
I_3:=\sum_{i\geq 1} \frac{2^{i-1}}{3^i} 
E \left[p_{\tT_i} \right], &&
I_4:=\sum_{i\geq 1} \frac{2^{i-1}}{3^i} 
E \left[p_{\tT_i-\tT_{i-1}} \right].
\end{eqnarray}
Since $\tT_i - \tT_{i-1}$ is exponentially
distributed with parameter $3$ (for all $i\geq 1$),
\begin{eqnarray}
I_4 &=& \sum_{i\geq 1} \frac{2^{i-1}}{3^i} \int_0^\infty ds 3e^{-3s} 
\frac{1-e^{-2s}}{2}
= \frac{3}{2} \left(\frac{1}{3} -\frac{1}{5} \right)=\frac{1}{5}.
\end{eqnarray}
Next, since $\tT_i$ follows a $\Gamma(i,3)$-distribution,
\begin{equation}
I_3 = \sum_{i\geq 1} \frac{2^{i-1}}{3^i} \int_0^\infty ds \frac{3^i}{(i-1)!} 
s^{i-1} e^{-3s} \frac{1-e^{-2s}}{2} 
=\int_0^\infty ds e^{-s} \frac{1-e^{-2s}}{2}= \frac{1}{3},
\end{equation}
and (using the substitution $u=e^{-s}$)
\begin{eqnarray}
I_1 &=& \sum_{i\geq 1} \frac{2^{i-1}}{3^i} \int_0^\infty ds \frac{3^i}{(i-1)!} 
s^{i-1} e^{-3s} \frac{1-e^{-2s}}{1+e^{-2s}} 
=\int_0^\infty ds e^{-s} \frac{1-e^{-2s}}{1+e^{-2s}}\ala
&=& \int_0^1 du \frac{1-u^2}{1+u^2}=2 \arctan 1 -1 = \frac{\pi}{2}-1.
\end{eqnarray}
Finally, using the independence between $\tT_{i-1}$ and $\tT_i-\tT_{i-1}$,
we get
\begin{eqnarray}
I_2 &=& \frac{1}{3}E\left[\frac{p_{\tT_1}}{1-p_{\tT_1}} \right] \ala
&&+\sum_{i\geq 2} \frac{2^{i-1}}{3^i} \int_0^\infty ds \frac{3^{i-1}}
{(i-2)!}  s^{i-2} e^{-3s} \int_0^\infty dt 3e^{-3t} 
\frac{1-e^{-2t}}{1+e^{-2s-2t}} \ala
&=& \frac{1}{3} \int_0^\infty ds 3e^{-3s} \frac{1- e^{-2s}}{1+e^{-2s}}
+\int_0^\infty ds\int_0^\infty dt 2 e^{-s}e^{-3t} 
\frac{1-e^{-2t}}{1+e^{-2s-2t}} \ala
&=& \int_0^1 du \frac{u^2(1-u^2)}{1+u^2} 
+ 2 \int_0^1 du \int_0^1 dv v^2 \frac{1-v^2}{1+u^2v^2}\ala
&=& \left(\frac{5}{3} - \frac{\pi}{2}\right) + 2 \left(\frac{\pi}{4}
- \frac{2}{3} \right)= \frac{1}{3}.
\end{eqnarray}
We finally get that $I = \frac{\pi}{2}-1 + \frac{1}{3}
+ \frac{1}{6} + \frac{1}{10}=\frac{\pi}{2}- \frac{2}{5}>1$. Thus
$E[Y_1] <0$.

\vip

{\bf Step 5.}
We still have to prove the exponential moment estimate. Using
the same notation as previously, we will just use that
for each $n \geq 1$,
$Y_1 \leq X_n=1+Z_n$ on $\Omega_n$. Recalling (\ref{jablat})
and that $\sum_{n\geq 1} P[\Omega_n]=1$, we get, for any
$k \geq 1$
\begin{equation}
P[Y_1  \geq k ] \leq \sum_{n\geq 1} P[Z_n \geq k-1, \Omega_n]
\leq 
\frac{1}{2^{k-1}}.
\end{equation}
We classically conclude that for $\e\in(0,\ln 2)$, $E[e^{\e Y_1}]<\infty$.

\vip
{\bf Step 6.} We finally have to build $\tzeta_0^1$. First note that
obviously, $\tzeta_{{\tT_1^1}}(i+\tR_{{\tT_1^1}}-1/2) \leq 1=\tzeta_0^1(i)$ 
for $i\leq 0$, while 
$\tzeta_{{\tT_1^1}}(1+\tR_{{\tT_1^1}}-1/2)=0=\tzeta_0^1(1)$
due to Lemma \ref{1firstpp}. Hence we just have to build 
$\tzeta_0^1(i)$ for $i\geq 2$. 
We write for simplicity $K_i=i+\tR_{{\tT_1^1}}-1/2$. 

\vip

Using the same arguments as in Step 3, one may check that conditionally
to $\cH_{\tT_1^1}$, the family 
$(\tzeta_{{\tT_1^1}}(K_i))_{i \geq 2}$ is independent,
and that for all $i \geq 2$,
\begin{equation}\label{neededlater}
P \left[ \left.\tzeta_{{\tT_1^1}}(K_i)=1 \right\vert \cH_{\tT_1^1} \right] 
\leq 1/2.
\end{equation}
We consider a family $(U_i)_{i \geq 2}$ of i.i.d. random variables uniformly
distributed on $[0,1]$ (independent of everything else) 
and we set, for each $i \geq 2$,
\begin{equation}\label{dfzetazero}
\tzeta_0^1(i):=\indiq_{\{\tzeta_{{\tT_1^1}}(K_i)=1\}} + 
\indiq_{\{\tzeta_{{\tT_1^1}}(K_i)=0\}} \indiq_{\{ U_i < \epsilon(i) \}},
\end{equation}
where, due to (\ref{neededlater}),
\begin{equation} 
\epsilon(i) = \frac{1 - 2 P[\tzeta_{{\tT_1^1}}(K_i)=1 \vert 
\cH_{{\tT_1^1}}]}{ 2 P[\tzeta_{{\tT_1^1}}(K_i)=0 \vert \cH_{{\tT_1^1}}]}
\in [0,1].
\end{equation}
We next observe that for all $i\geq 2$, 
\begin{eqnarray}
P\left[\tzeta_{0}^1(i)=1 \vert \cH_{{\tT_1^1}}\right] 
 =  E\left[P[\tzeta_{0}^1(i)=1 \vert \cH_{{\tT_1^1}}] \right]
\hskip3cm\ala
= E\left[ 
P[\tzeta_{{\tT_1^1}}(K_i)=1 \vert \cH_{{\tT_1^1}}]   
+ \e (i) P[\tzeta_{{\tT_1^1}}(K_i)=0 \vert \cH_{{\tT_1^1}}]\right]=\frac{1}{2}
\end{eqnarray}
due to our choice for $\e(i)$ and to the independence of $U_i$ of
everything else. We deduce that for each $i\geq 2$,
$\tzeta_{0}^1(i)$ is a Bernoulli random variable with parameter
$1/2$, and that it is independent of $\cH_{{\tT_1^1}}$.
This and the conditionnal (to $\cH_{{\tT_1^1}}$) independence of
the family $(\tzeta_{0}^1(i))_{i \geq 2}$ clearly imply that
finally, $(\tzeta_{0}^1(i))_{i\geq 2}$ is an i.i.d.
sequence of Bernoulli random variables with parameter $1/2$, independent
of $\cH_{{\tT_1^1}}$. Finally, it is clear from (\ref{dfzetazero})
that for all $i\geq 2$, $\tzeta_0^1(i) \geq \tzeta_{{\tT_1^1}}(K_i)$.
This concludes the proof.
\end{proof}

The following lemma shows a way to bound from above
the right contour process started with the initial 
condition $\zeta_0\sim\Gamma$
by a continuous-time random walk.

\begin{lem}\label{upperrw}
Let $\zeta_0 \sim \Gamma$, let $N,V$ be two independent IFPPs, and consider
the $(\zeta_0,N,V)$-right contour process $(R^0_t)_{t\geq 0}$ around $0$.
Then for $k \geq 0$, $P[R^0_0=k-1/2]=(1/2)^{k+1}$. Furthermore, we may find
a Poisson process $(Z_t)_{t\geq 0}$ with rate $1$, a family
of i.i.d. random variables $(Y_i)_{i\geq 1}$ distributed as $Y_1$
(see Lemma \ref{renew}) in such a way that $R_0^0$ and $((Z_t)_{t\geq 0},
(Y_i)_{i\geq 1})$ are independent, while a.s., for all $t\geq 0$,
\begin{equation}\label{ineqcool}
R^0_t \leq R^0_0+ \sum_{i=1}^{Z_t} Y_i.
\end{equation}
\end{lem}

\begin{proof} 
We omit as in the proof of Lemma \ref{renew} the superscript $0$.
We consider $\zeta_0\sim \Gamma$ to be fixed, and a $\zeta_0$-right contour
process $(R_t)_{t\geq 0}$ around $0$. First,
it is obvious that for $k \geq 0$,
\begin{equation}
P[R_0=k-1/2]=P[\zeta_0(0)=1,...,\zeta_0(k-1)=0,\zeta_0(k)=1]=(1/2)^{k+1},
\end{equation}
since $\zeta_0 \sim \Gamma$. Next, let us explain (\ref{ineqcool}).
The main ideas are the following: we first bound our contour
process by a contour process $\tR^1$ which starts from a 
(shifted) $\Xi$-distributed
initial data. When this process first jumps to the right,
at some instant $\tau_1$, we bound the Bernoulli process at this time
by a shifted $\Xi$-distributed data, independent of 
$(\tR^1_t)_{t\in[0,\tau_1]}$. Thus we make start again a contour process
$\tR^2$ from this $\Xi$-distributed initial data in such a way
that it dominates (with a shift) $\tR^1$, and thus $R$. And so on...
The advantage of this method is that the 
{\it increments} (between two {\it renewal times}) are independent.

\vip

We define $\zeta_0^1$ by
$\zeta_0^1(i)=1$ if $i \leq 0$, $\zeta_0^1(1)=1$,
and $\zeta_0^1(i)=\zeta_0(i+R_0-1/2)$ for $i \geq 2$. We observe that 
$\zeta_0(i+R_0-1/2)\leq \zeta_0^1(i)$ for all $i\in \zz$, 
using Lemma \ref{1firstpp}. We also notice that $\zeta_0^1$ is
independent of $R_0$, and is $\Xi$-distributed.

We thus may find, using Lemma \ref{toctoc}, a contour process
$(\tR^1_t)_{t\geq 0}$ around $0$, 
starting from $\zeta_0^1$, independent of $R_0$,
such that for all times, $R_t- R_0 \leq \tR^1_t-1/2$
(recall that $R_t-R_0$ starts from $1/2$, and that for all $i\in\zz$,
$\zeta_0(i+R_0-1/2)\leq\zeta^1_0(i)$). On the other hand, we consider
$\tau_1=\inf \{t\geq 0; \Delta \tR^1_t>0 \}$ the first instant
where $\tR^1$ jumps to the right, so that $\tR^1_t \leq \tR^1_{\tau_1}$ 
for all $t\in[0, \tau_1]$.
Hence setting $Y_1:=\tR^1_{\tau_1}-1/2$, we finally
obtain that a.s., for all $t\in[0, \tau_1]$, $R_t \leq R_0 + Y_1$.
We also observe that $R_0$ and $(\tau_1,Y_1)$ are independent.
Finally, $\tau_1$ is  exponentially distributed with parameter $1$,
due to Lemma \ref{firstpp}-(d), and $Y_1$ is distributed as in Lemma
\ref{renew} by construction.

\vip

Due to Lemma \ref{renew}, we may find $\zeta^2_0 \sim \Xi$,
independent of $R_0$ and $Y_1$, such that a.s., for all $i \in \zz$,
$\zeta^1_{\tau_1}(i+\tR^1_{\tau_1}-1/2) \leq \zeta^2_0(i)$,
where $(\zeta^1_t)_{t \geq 0}$ is the Bernoulli process
starting from $\zeta^1_0$ associated with the contour process 
$(\tR^1_t)_{t\geq 0}$.

Using Lemma \ref{toctoc}, we thus may build a contour process
$(\tR^2_t)_{t\geq 0}$ around $0$, independent of $R_0$ and $Y_1$, such that
for all times $\tR^1_{\tau_1+t} - \tR^1_{\tau_1} \leq \tR^2_t -1/2$.
As a consequence, we observe that for all $t\geq 0$,
\begin{equation}
R_{\tau_1+t} -R_0 \leq \tR^1_{\tau_1+t}-1/2 \leq (\tR^1_{\tau_1}-1/2)
+ (\tR^2_t -1/2) \leq Y_1 +  (\tR^2_t -1/2).
\end{equation}
Denote by $\tau_2=\inf \{t\geq 0; \Delta \tR^2_t>0 \}$ the first instant
where $\tR^2$ jumps to the right, so that $\tR^2_t \leq \tR^2_{\tau_2}$ 
for all $t\in[0, \tau_2]$.
Hence setting $Y_2:=\tR^2_{\tau_2}-1/2$, we finally
obtain that a.s., for all $t\in[\tau_1,\tau_1+\tau_2]$, $R_t \leq R_0 + Y_1
+Y_2$.
We also observe that $(\tau_2,Y_2)$ is independent of $R_0$ and 
$(\tau_1,Y_1)$.
Finally, $\tau_2$ is  exponentially distributed with parameter $1$,
due to Lemma \ref{firstpp}-(d), and $Y_2$ is distributed as $Y_1$
by construction.

\vip

Iterating the procedure, we find an i.i.d. family $(\tau_k,Y_k)_{k\geq 1}$
of random variables, independent of $R_0$,
such that $Y_1$ is distributed as in Lemma \ref{renew} and
$\tau_1$ is exponentially distributed with parameter $1$,
such that for all $t\geq 0$, 
\begin{equation}
R_t - R_0  \leq \sum_{k\geq 1} Y_k \indiq_{\{t \geq \tau_1+...+\tau_k\}}.  
\end{equation}
This ends the proof.
\end{proof}

We finally conclude the proof of the main estimates of this section.

\vip

\begin{proof} {\bf of Proposition \ref{est}.}
We omit the superscript $0$ for simplicity.
The proof is based on the use of Lemmas \ref{upperrw} and \ref{renew}.
We thus write, according to (\ref{ineqcool}),
$R_t \leq R_0 + S_{Z_t}$, with $S_n=Y_1+...+Y_n$.

First of all, we deduce from Lemma \ref{renew}-(i)-(ii)
that there exists $\gamma\in (0,\ln 2)$ 
such that $q :=E[e^{\gamma Y_1}]<1$.

Since $\gamma \in (0,\ln 2)$, we also deduce from Lemma \ref{upperrw}
that $E[e^{\gamma R_0}]<\infty$.

Next, we recall that since $(Z_t)_{t\geq 0}$ is a Poisson process
with rate $1$, 
for all $t\geq 0$, $P[Z_t \leq t/2] \leq e^{-\delta t}$, where
$\delta:=(1-\ln 2)/2>0$ (any $\delta>0$ would work as well). 


We have
$\oR_\infty=\sup_{t \geq 0} R_t \leq R_0 + \sup_{n \geq 1} S_n$.
This implies that $e^{\gamma \oR_\infty} 
\leq e^{\gamma R_0} \sum_{n \geq 1} 
e^{\gamma S_n}$.
Thus, since $q\in (0,1)$, and since $S_n$ and $R_0$ are independent, 
\begin{equation}
E[e^{\gamma \oR_\infty}] \leq E[e^{\gamma R_0}]
\sum_{n \geq 1} E[e^{\gamma S_n}] \leq C \sum_{n\geq 1} E[e^{\gamma Y_1}]^n
= C \sum_{n \geq 1} q^n <\infty.
\end{equation}
Next we want to upperbound $\rho$. First, 
$\rho \leq \inf\{t\geq 0, R_t <0$ and $L_t>0\}$, so that
by symmetry, for any $t \geq 0$,
\begin{equation}
P[\rho \geq t] \leq P[R_t> 0 \hbox{ or } L_t< 0] \leq
2  P[R_t> 0].
\end{equation}
But, since $R_t\leq R_0+S_{Z_t}$,
\begin{eqnarray}
&&P[\rho \geq t] \leq 2 P[R_t > 0] \leq 2P[ Z_t \leq t/2] + 
2P[R_0 + \sup_{n > t/2} S_n >0  ] \ala
&&\leq 2e^{-\delta t} + 2 E \left[e^{\gamma(R_0 + \sup_{n > t/2} S_n)} \right]
\leq  2 e^{-\delta t}+ 2 E\left[e^{\gamma R_0}\right] \sum_{n> t/2}E\left[
e^{\gamma S_n} \right]\ala
&&\leq   2 e^{-\delta t} + C \sum_{n>t/2} q^n \leq 2e^{-\delta t} +
C q^{t/2} \leq Ae^{-at}
\end{eqnarray}
for some constants $A>0$, $a>0$.
We classically conclude that
for any $\beta \in (0,a)$, $E[e^{\beta \rho}]<\infty$.
\end{proof}

\section{Proof of Theorem \ref{main}}\label{concl}\setcounter{equation}{0}

Our aim in this section is to conclude the proof of our main result.

\vip

\begin{proof} {\bf of Theorem \ref{main}}
We divide the proof into several steps. 
Let us recall briefly the notation of Proposition \ref{ps1} and
Definition \ref{dfcp}: we consider two independent IFPPs 
$N,V$, and $\zeta_0 \sim \Gamma$.  
Let $(\zeta_t)_{t\geq 0}$ be the 
$(\zeta_0,N)$-Bernoulli process $(\zeta_t)_{t\geq 0}$, and let
$(\tzeta_t)_{t\in(-\infty,0]}$ be its time-reversed
built in Lemma \ref{revtimeber}.

For $T \in(-\infty,0]$ and $\varphi \in E_{\zeta_T}$ (recall (\ref{dfez})), 
we denote by
$(\tzeta_t,\eta^{T,\varphi}_t)_{t\in [T,0]}$ the 
$(\tzeta_T,\varphi,N^T,V^T)$-coupled 
Bernoulli avalanche process 
with $N^T_t(i)=N_{(-T)-}(i)-N_{(-t)-}(i)$ and 
$V^T_t(i)=V_{(-T)-}(i)-V_{(-t)-}(i)$
for $t\in [T,0]$ and $i \in \zz$. 
Recall that 
\begin{equation}
\tau_i=\sup \{T \leq 0 ;
\forall \varphi\in E_{\zeta_T},\; \eta^{T,\varphi}_0(i)=\eta^{T,\bZ}_0(i)\}.
\end{equation}
We then may put $\eta_0(i):=\eta^{\tau_i+s,\bZ}_0(i)$ (for some $s<0$,
recall Proposition \ref{ps1})
provided $\tau_i>-\infty$.
Recall also
that by definition of $\tau_i$, we have, for all $T<0$,
all $\varphi \in E_{\zeta_T}$, $\eta^{T,\varphi}_0(i)=\eta_0(i)$
on the event $\{T<\tau_i\}$.

Next, we consider
the $(\zeta_0,N,V)$-left and right contour processes $(L^i_t)_{t\geq 0}$
and $(R^i_t)_{t\geq 0}$ around $i$, for each $i \in \zz$, and we
adopt the notation 
\begin{eqnarray}\label{barbar}
\oR^i_t=\sup_{s \in [0,t]} R^i_s, \quad
\uL^i_t = \inf_{s \in [0,t]} L^i_s,\quad
\rho^i = \inf\{t \geq 0;\; R^i_t <  L^i_t \}.
\end{eqnarray}
Finally, 
for $k<l \in \zz\cup\{-\infty,+\infty\}$ and $t\in [0,\infty]$, 
we consider the $\sigma$-field
\begin{equation}
 \cG_{k,l,t} =\sigma \left\{\zeta_0(j), \; V_s(j),\; N_s(j),\; s\in [0,t],\; 
l\leq j \leq k  \right\}. 
\end{equation}

{\bf Step 1.} We first show that a.s., $-\tau_i \leq \rho^i$ a.s.,
and that $\eta_0(i)= \Phi_i (Z^i)$, for some deterministic function
$\Phi_i$, where 
\begin{equation}
Z^i:=\left(\zeta_0(j), N_s(j),V_s(j), j\in
\{\uL^i_{\rho^i}-1/2,...,\oR^i_{\rho^i}+1/2\}, s\in [0,\rho^i]\right).
\end{equation}
The function $\Phi^i$ can not easily be made explicit,
see however Step 2 of the algorithm described in
the Appendix \ref{appendix}.

It clearly suffices to treat the case $i=0$. 
We consider the {\it box} delimited by $L^0_t$ and $R^0_t$ until they meet,
i.e. until $t=\rho^0$. Consider an avalanche process
starting at some time $T< - \rho^0$ with a given initial condition
$\varphi \in E_{\tzeta_T}$. We wish to rebuild its value at time $0$,
thus the time goes now down on Figure 2 (see also Figure 4 below).

Observe, having a look at Figure 2, 
that on the right and left of this box,
the Bernoulli process $\tzeta_t$ is vacant (see Lemma \ref{1firstpp}),
so that due to our coupling, the avalanche process is also vacant,
since it is always smaller (see Proposition \ref{avaber}). 
Hence no interaction can go inside this box
(from its left and right sides),
since vacant sites {\it cut} the interaction: indeed, flocks
falling outside this box can not make die flocks inside the box.

Next, notice that
the horizontal segments delimiting this box on the top side contain
sites of the following type: 

(a) either vacant sites of the Bernoulli process,
so that the avalanche process is also vacant at this site 
at this time (since it is always smaller);

(b) either sites where the Bernoulli process becomes occupied because 
of a black mark, so that the avalanche process also becomes occupied 
at this time (because when the Bernoulli process is vacant, then the avalanche
and Bernoulli processes both become occupied  when they encounter
a black mark);

(c) either a grey mark (in the middle of an occupied zone of the Bernoulli
process), so that at this site and at this time, the avalanche process is
(or becomes) occupied.

As a conclusion, the avalanche process is vacant on the left and 
right of the box, and the value
of the avalanche process on the top horizontal segments of
this box are determined, independently of its starting time
$T<-\rho^0$ and initial condition $\varphi \in E_{\tzeta_T}$.
We thus may rebuild the avalanche process 
$(\eta_t^{T,\varphi})_{t \in [\-\rho^0,0]}$, and the obtained
value $\eta^{T,\varphi}_0(0)$ does not depend on $T<-\rho^0$ nor on
$\varphi \in E_{\tzeta_T}$. We thus deduce that 
$\tau_0 \geq -\rho^0$ and that $\eta_0(0)=\eta^{T,\varphi}_0(0)$.
This value $\eta_0(0)$
clearly depends only on the values of $N,V,\tzeta_t$ in this box,
and the fact that the Bernoulli process is vacant on the outside boundary of
this box. We thus can say that $\eta_0(0)$ is a (deterministic) function
of $Z^0$.

\begin{figure}[t]\label{figrecon}
\centerline{Figure 4: Reconstruction of $\eta_0$ (Step 1).}
\vskip0.1cm
\centerline{ \hskip1cm \includegraphics{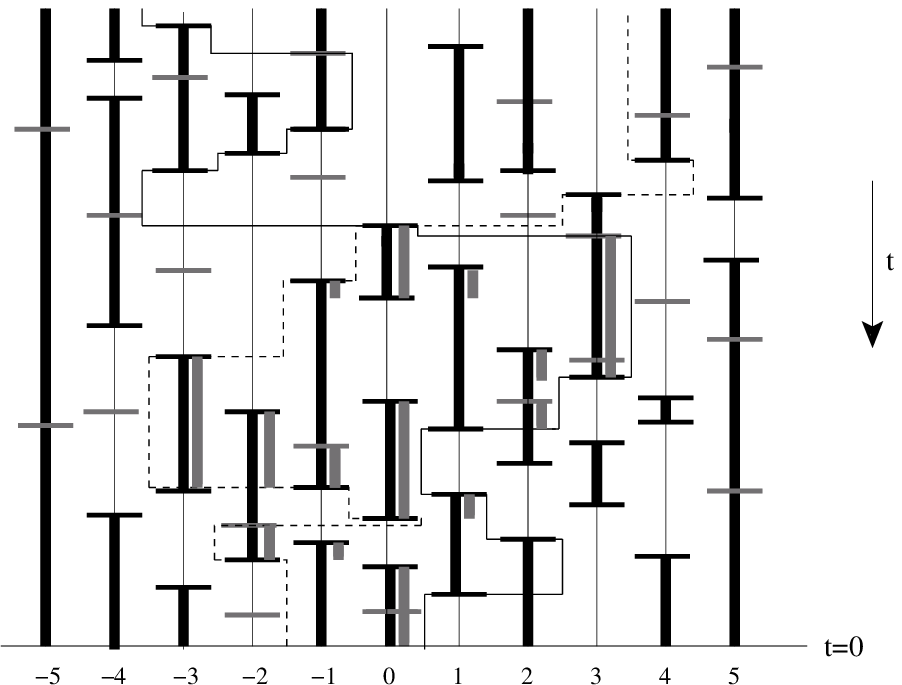} }
\begin{footnotesize}
We explain here graphically how to obtain the values of 
$\eta_0(i)$ for $i \in \{-2,..,1\}$. The avalanche process
is represented on the right of each site, in grey.

We start from the top of the box, that is at site $0$.
A flock appears here, so that $0$ becomes occupied (independently of
its starting time $T<-\rho^0$ and initial condition 
$\varphi \in E_{\tzeta_T}$). 

Next, a grey mark appears at the site $3$, 
which thus remains occupied (if it was
already) or becomes occupied (if it was not), so that in any case,
the site $3$ is occupied at this time.

At the same time,
the sites $1$ and $2$ are vacant, since they are vacant for the 
Bernoulli process.

Next, $1$ and then 
$-1$ become occupied due to black marks. 
But then the flock at $0$ dies, and makes 
$-1,0,1$ become vacant. And so on... 

This way, we see that the
site $-1$ is finally vacant at time $0$, while the site $0$ is finally 
occupied. On the other hand, it is immediate that $-2$ and $1$
are vacant, since the avalanche process is smaller than the Bernoulli
process with our coupling.

We could also see that $\eta_0(2)=1$ and $\eta_0(3)=0$ here, but
it is not possible to decide if $\eta_0(-3)=1$, because
it could be killed by a flock dying at some site $i\leq -6$.
\end{footnotesize}
\end{figure}

\vip

{\bf Step 2.} We also observe that for $i\in\zz$, $k\leq i \leq l$, 
for $T>0$,  
\begin{equation}
\Omega_i(k,l,T):=\left\{\rho^i<T,  \uL^i_{T} \geq k+1/2, 
\oR^i_{T}\leq l-1/2 \right\}
\in \cG_{k,l,T}
\end{equation}
and that $Z^i \indiq_{\{\Omega_i(k,l,T)\}}$ is $\cG_{k,l,T}$-measurable.
This is clear from Figure 2.

\vip

{\bf Step 3.} Next, we notice that point (a) (existence and uniqueness
of an invariant distribution $\Pi$ for the avalanche process) follows 
immediately from 
Proposition \ref{ps1}, provided we know that $\tau_i>-\infty$ a.s.~for 
all $i\in\zz$. But we know from Step 1 that $\tau_i \geq -\rho^i$.
Lemma \ref{firstpp}-(a) implies that for all $i \in \zz$, 
$\rho^i$ and $\rho^0$ are identically distributed. As a consequence,
it suffices to show that $\rho^0<\infty$ a.s., which follows
from Proposition \ref{est}-(a).

Hence, $\Pi=\cL(\eta_0)$ is the unique invariant distribution of the 
avalanche process.

\vip

{\bf Step 4.} Point (c) (existence of a perfect simulation algorithm
for $(\eta_0(i))_{i=-l,...,l}$) is also immediate, see Appendix for an explicit
simulation algorithm. Let $l\geq 0$ be fixed. 
Due to Step 3, we know that $\eta_0\sim \Pi$, it thus suffices to  
simulate (perfectly) 
$(\eta_0(j))_{j \in \{-l,...,l\}}$. This can be done by simulating 
first $Z^{-l},...,Z^l$.
This can be done due to Step 2, which says that $Z^i$ depends on
$\zeta_0$, $N$, and $V$ on an a.s.~finite number of sites and on an 
a.s.~finite time intervall. Next, it suffices to compute 
$\eta_0(i)=\Phi_i(Z^i)$ for all $i \in \{-l,...,l\}$, which can be done
following the rules explained in Step 1, see also Figure 4.

\vip

{\bf Step 5.} We now check the mixing property anounced in point (d).
Let thus $n\geq 1$ and $k\in \zz$ be fixed. 
We consider the events (here $[x]$ stands for the integer part of
$x\in\rr$)
\begin{equation}
\Omega^1_n=\left\{\oR^k_\infty \leq k+[n/3]-1/2\right\} \hbox{ and } 
\Omega^2_n=\left\{
\uL^{k+n}_\infty \geq k+[2n/3]+1/2 \right\}.
\end{equation}
We deduce from Lemma \ref{firstpp}-(c) that on 
$\Omega^1_n$, $\oR^i_\infty\leq k+[n/3]-1/2$
for all $i \leq k$.  We thus deduce from Step 2 that the family
$(Z^i \indiq_{\Omega^1_n})_{i\leq k}$ is 
$\cG_{-\infty,k+[n/3],\infty}$-measurable.
Hence $(\eta_0(i) \indiq_{\Omega^1_n})_{i\leq k}$ is 
$\cG_{-\infty,k+[n/3],\infty}$-measurable.
By the same way, $(\eta_0(i) \indiq_{\Omega^2_n})_{i\geq k+n}$ 
is $\cG_{k+[2n/3],\infty,\infty}$-measurable. 

But of course, the two $\sigma$-fields $\cG_{-\infty,k+[n/3],\infty}$ 
and $\cG_{k+[2n/3],\infty,\infty}$ are independent. Hence,
for any $A \subset \{0,1\}^{(-\infty,k]}$, $B \subset \{0,1\}^{[k+n,\infty]}$, 
\begin{eqnarray}
\Big|P\left[(\eta_0(i))_{i \leq k}\in A, 
(\eta_0(i))_{i \geq k+n}\in B \right]\hskip3cm \ala 
- P\left[(\eta_0(i))_{i \leq k}\in A\right]
P\left[(\eta_0(i))_{i \geq k+n}\in B \right]\Big| \ala
\leq 2P[(\Omega_n^1)^c]+ 2P[(\Omega_n^2)^c],
\end{eqnarray}
so that
\begin{eqnarray}
\left|\Pi_{(-\infty,k]\cup [k+n,\infty)} -
\Pi_{(-\infty,k]} \otimes \Pi_{[k+n,\infty)}\right|
\leq 2P((\Omega_n^1)^c)+ 2P((\Omega_n^2)^c).
\end{eqnarray}
Using Lemma \ref{firstpp}-(a), we deduce that $P[(\Omega_n^1)^c]
=P[\oR^k_\infty \geq k+[n/3]+1/2]=P[\oR^0_\infty \geq [n/3]+1/2]$ and
$P[(\Omega_n^2)^c]=P[\uL^{k+n}_\infty \leq k+[2n/3]-1/2]=P[\oR^0_\infty \geq 
n-[2n/3]+1/2]\leq P[\oR^0_\infty \geq [n/3]+1/2]$, since 
$n-[2n/3]\geq [n/3]$. Using finally
Proposition \ref{est}-(b), obtain that for some constants $\gamma>0$, $K>0$,
\begin{eqnarray}
P[(\Omega_n^1)^c] + [(\Omega_n^2)^c] &\leq& 2 P[\oR^0_\infty 
\geq [n/3]+1/2]\ala
&\leq& 2 e^{-\gamma/2-\gamma [n/3]} E[e^{\gamma \oR^0_\infty}] 
\leq K e^{-\gamma n /3}.
\end{eqnarray}
We thus obtain (\ref{mix}), setting $q=e^{-\gamma/3}$ and $C=2K$.
 
\vip

{\bf Step 6.} Finally, it remains to study the trend to equilibrium.
For $\varphi \in E$, $t\geq 0$, we denote by $\Pi^\varphi_t$
the law of the $\varphi$-avalanche process at time $t$.
The main difficulty here is to obtain the trend to equilibrium for any initial
datum $\varphi \in E$, since our coupling allows us {\it a priori}
to deal only with initial data stochastically smaller than $\Gamma$.
We thus have to introduce a final coupling, which mixes those
of Lemma \ref{bermono} and Proposition \ref{ps1}. We keep, however,
all the notations introduced in this proof.
We fix $\varphi \in E$.

{\it Step 6.1.} We consider a third IFPP $W$, independent of $N$, $V$,
and $\zeta_0$. We consider, for $T<0$ using Lemma \ref{bermono},
the $(\tzeta_T,\varphi,N^T,W^T)$-coupled Bernoulli processes
$(\tzeta_t,\bzeta^T_t)_{t\in [T,0]}$.  Recall that
$(\tzeta_t)_{t\in [T,0]}$ is the $(\varphi, N^T)$-Bernoulli process,
while
$(\bzeta_t^T)_{t\in [T,0]}$ is the $(\varphi, M_T)$-Bernoulli process
for some IFPP $M_T$. Notice that due to Lemma \ref{bermono}-(ii)-(iii),
we obtain that for any $k\geq 0$, any $t\in [T,0]$,
\begin{eqnarray}\label{xxx1}
&\hbox{if } \Omega_{1}(t,T,k):=
\{\bzeta_s^T(i)=\tzeta_s(i);
\forall \; i\in[-k,k], s \in [t,0]\} \ala
&\hbox{then } P[\Omega_{1}(t,T,k)] \geq 1 - (2k+1)e^{-2(t-T)}.
\end{eqnarray}
We finally consider the $(\varphi,\varphi,M_T,V^T)$-coupled Bernoulli-avalanche
processes $(\bzeta_t,\bareta_t^T)_{t\in [T,0]}$.

{\it Step 6.2.} 
We consider the event, for $n \geq 0$,  $t>0$,
\begin{equation}
\Omega_0(n,t)=\left\{\oR^l_\infty \leq l+n+1/2, \uL^{-l}_\infty \geq -l-n-1/2, 
\max_{i\in\{-l,...,l\}}\rho^i < t \right\}.
\end{equation}
We know that on this event, $\eta_0^{-t,\varphi}(i)=\eta_0(i)$
for all $i \in \{-l,...,l\}$, as soon as $\varphi \in E_{\zeta_{t}}$.
But we easily understand, using Step 6.1 and some arguments as
in Step 1, that on $\Omega_{1}(-t,-2t,l+n+1)\cap \Omega_0(n,t)$,
we also have $\bareta_0^{-2t}(i)=\eta_0(i)$ for all $i \in \{-l,...,l\}$.

{\it Step 6.3.}
On the other hand, $\cL(\bar\eta_0^{-2t})=\Pi^\varphi_{2t}$ and
$\cL(\eta_0)=\Pi$, so that we classically 
deduce that for all $n\geq 0$,
\begin{equation}\label{xxx2}
\left|(\Pi^\varphi_{2t})_{[-l,l]}-\Pi_{[-l,l]} 
\right|_{TV} \leq 2 P[(\Omega_0(n,t))^c\cup (\Omega_{1}(-t,-2t,l+n+1))^c ].
\end{equation}

{\it Step 6.4.}
We obtain, using Proposition \ref{est} and Lemma 
\ref{firstpp}-(a), that
\begin{eqnarray}\label{xxx3}
P[(\Omega_0(n,t))^c ]&\leq& P \left[\oR^l_\infty \geq l+n+3/2 \right]
+ P \left[ \uL^{-l}_\infty \leq -l-n-3/2 \right] \ala
&& + \sum_{i=-l}^l P[\rho^i\geq t] \ala
&\leq & 2 \left[\oR^0_\infty \geq n+3/2 \right]+ (2l+1)P[\rho^0\geq t]\ala
&\leq & 2e^{-\gamma n} E\left[e^{\gamma \oR^0_\infty} \right] 
+ (2l+1) e^{-\beta t} E\left[ e^{\beta \rho^0} \right] \ala
&\leq&
A(e^{-\gamma n} + (2l+1)e^{-\beta t}).
\end{eqnarray}
for some constants $A>0$ $\beta>0$, $\gamma>0$.

{\it Step 6.5.} Gathering (\ref{xxx1}), (\ref{xxx2}), (\ref{xxx3}),
we finally obtain
that for any $t\geq 0$, any $n\geq 1$,
\begin{eqnarray}
\left|(\Pi^\varphi_{2t})_{[-l,l]}-\Pi_{[-l,l]} \right|_{TV} 
&\leq& 2A(e^{-\gamma n} + (2l+1)e^{-\beta t}) \ala
&&+ 2(2(l+n+1)+1)e^{-2t}.\ala
\end{eqnarray}
Choosing finally $n=[t]$, we deduce that for $a=\min(1,\gamma,\beta)$,
we may find a constant $K>0$ such that for all $t>0$,
all $l \geq 0$,
\begin{equation}
\left|(\Pi^\varphi_{2t})_{[-l,l]}-\Pi_{[-l,l]}\right|_{TV} \leq K(1+l) 
e^{-a t}.
\end{equation}
This yields (\ref{tteap}), and concludes the proof.
\end{proof}

\section{A related mean-field model}
\label{cf} \setcounter{equation}{0}

This section, quite independent of the rest of the paper, is devoted
to the brief study of a mean-field coagulation-fragmentation model
related to the avalanche process.

To obtain a process which preserves the total mass, we will slightly
change our point of view: we assume that each edge of $\zz$
has a mass equal to $1$. 

Consider a possible state $\eta \in E$ of the avalanche process.
Two neighbour edges, say $(i-1,i)$ 
and $(i,i+1)$, are said to belong to the same {\it particle}
if $\eta(i)=1$: the flock lying at $i$ {\it glues} the two edges.
For example, the edge $(0,1)$ belongs to a particle with
mass $3$ if $\eta(-1)=\eta(2)=0$ and $\eta(0)=\eta(1)=1$,
or if $\eta(0)=\eta(3)=0$ and $\eta(1)=\eta(2)=1$.
Similarly, $(0,1)$ belongs to a particle with mass $1$ if and only
if $\eta(0)=\eta(1)=0$.

\vip

Then we consider, for $\eta \in E$ and for $k\in \nn$, if it exists,
\begin{equation}
c_k(\eta)= \lim_{n \to \infty} \frac{ 
\hbox{number of particles with mass } k \hbox{ in } [-n,n] }{2n+1},
\end{equation}
which represents the average number of particles with mass $k$
per unit of length.
Consider an avalanche process $(\eta_t)_{t\geq 0}$.
Assume for a moment 
that for each $t\geq 0$, the successive masses of particles in
$\eta_t$ are independent
(which is intuitively far from being exact).
Then, using the invariance by translation of the model,
one would have, for $k\geq 1$, $t\geq 0$,
\begin{equation}
c_k(t):=c_k(\eta_t)=\frac{1}{k}P \left[(0,1) \hbox{ belongs to a 
particle with mass } k \hbox{ in } \eta_t \right]
\end{equation}
The family $(c(t))_{t\geq 0}=(c_k(t))_{t\geq 0,k\geq 1} $
would also satisfy   $\sum_{k \geq 1} k c_k(t)=1$ for all $t\geq 0$, and
the following infinite system of differential equations:
\begin{eqnarray}\label{coagfrag}
\frac{d}{dt} c_1(t) &=&- 2 c_1(t)+ \sum_{k\geq 1} (k-1)kc_k(t) ,\\
\frac{d}{dt} c_k(t) &=&- 2 c_k(t)- (k-1)c_k(t) + \frac{1}{m_0(c(t))} 
\sum_{i=1}^{k-1}c_i(t)c_{k-i}(t) \; \; (k\geq 2),\nonumber
\end{eqnarray}
where $m_0(c(t))=\sum_{k \geq 1} c_k(t)$ is the average number of particles
per unit of length. Indeed, 
the first equation expresses that 
an isolated edge merges with its two neighbours with rate $1$,
while each time a flock falls on a particle with mass $k$,
which happens at rate $k-1$ (since a particle with mass $k$ contains
$k$ edges and thus $k-1$ sites), an avalanche occurs and 
$k$ new particles with mass $1$ appear.
Next, the second equation expresses that particles with mass $k$ become
larger at rate $2$ (when a flock falls on one of its two extremities),
that particles with mass $k$ disappear when they are subject to an avalanche
(which happens at rate $k-1$), and that particles with mass $k$ do
appear when a flock falls between a particle with mass $i$ and
a particle with mass $k-i$. This last event occurs at rate
$1$, proportionnaly to the number (per unit of length) 
of pairs of neighbour particles
with masses $i$ and $k-i$, which is exactly
$c_i(t) c_{k-i}(t)/m_0(c(t))$. We use the abusive independence assumption
when computing this last rate.

\vip

We refer to Aldous \cite[Construction 5]{aldous} 
for very similar considerations,
without fragmentation, where the independence between neighbours 
really holds.

\vip

The system (\ref{coagfrag}) can be seen as a coagulation-fragmentation
equation with constant coagulation rate $K(i,j)=2$, 
with a splitting rate (from one particle with mass $k$ into
$k$ particles with mass $1$) $F(k;1,...,1)=k-1$, 
the change of time $1/m_0(c(t))$ lying in front of the coagulation term.
Indeed, we could write, for example when $k \geq 2$,
\begin{eqnarray}
\frac{d}{dt} c_k(t) &=& - F(k;1,...,1) c_k(t) \\
&& + \frac{1}{m_0(c(t))}
\left[ - c_k(t) \sum_{i \geq 1} K(k,i)c_i(t) + 
\sum_{i=1}^{k-1}K(i,k-i) c_i(t)c_{k-i}(t) \right]. \nonumber
\end{eqnarray}
The term in brackets on the second line is the right-hand side member 
of the well-known Smoluchowski coagulation equation. 
See Aldous \cite{aldous}, Lauren\c cot-Mischler \cite{lm} 
for reviews on these types of equation. No result about trend
to equilibrium for such a model without detailed balance condition
(here the coagulation is binary, which is not the case
of fragmentation)
seem to be available. See however 
Fournier-Mischler \cite{fm} for partial 
results about a coagulation-fragmentation without balance condition
in the case of binary fragmentation.

\vip

We aim here to compute the steady state of this mean-field model,
and to show numerically that it approximates closely the invariant 
distribution of the avalanche process.

\begin{prop}\label{sscf}
The system of equations (\ref{coagfrag}) admits a unique
steady state $c=(c_k)_{k\geq 1}$, in the sense that: $c_k \geq 0$
for all $k \geq 1$, $\sum_{k \geq 1} k c_k =1$, and,
with $m_0(c)=\sum_{k \geq 1} c_k$, 
\begin{eqnarray}
2 c_1 &=& \sum_{k\geq 1} (k-1)kc_k,\ala
(k+1)c_k &=& \frac{1}{m_0(c)} 
\sum_{i=1}^{k-1}c_ic_{k-i} \; \; (k\geq 2).
\end{eqnarray}
This steady state is given by $c_k= a_k g^{k-1} 2^{-k}$,
where 
$a_1=1$ and for $k \geq 2$, $a_k=\frac{1}{k+1}\sum_{i=1}^{k-1} a_j a_{k-j}$,
while $g>0$ is the unique solution of
$\sum_{k\geq 1} a_k (g/2)^k =1$.

We also have $c_1= 1/2$, $\sum_{k \geq 1} k^2 c_k =2$,
and $m_0(c)=1/g$. 
\end{prop}

\begin{proof}
Consider the sequence $(a_k)_{k \geq 1}$ defined in the statement. Remark
that for any $x>0$, any $g>0$, the sequence $x_1=x$,
$x_k=\frac{g}{k+1} \sum_{i=1}^{k-1} x_j x_{k-j}$ is explicitely
given by $x_k=a_k x^k g^{k-1}$.

Thus, $(c_k)_{k\geq 1}$ is a steady state of (\ref{coagfrag})
if and only if there exist $x>0$ and $g>0$ such that

(i) $\forall$ $k\geq 1$, $c_k=a_k x^k g^{k-1}$,

(ii) $g=1/m_0(c)$, 

(iii) $\sum_{k \geq 1} k c_k =1$,

(iv) $ x= \frac{1}{2}
(\sum_{k \geq 1} k^2 c_k -1)$.

\vip

Points (i) and (ii) imply that necessarily, 
$\sum_{k\geq 1} a_k (xg)^k =1$. Thus $q:=xg$ is clearly uniquely defined,
and satisfies $0<q<1$ (since $a_1 =1$ and $a_2=1/3>0$).
Next, using (iii), we deduce that $g=\sum_{k\geq 1} k a_k q^k$ is
also uniquely defined (and finite, since $q<1$ and since one
easily checks recursively that $a_k \leq 1$ for all $k \geq 1$).
Thus $x=q/g$ is also uniquely defined. This shows that there exists at
most one steady state.
We next have to verify that these values for $x$ and $g$ imply point (iv).
Using the definition of $(a_k)_{k\geq 1}$, we obtain that
on the one hand,
\begin{equation}
\sum_{k \geq 2} (k+1)a_k q^k = \sum_{k \geq 2} q^k \sum_{j=1}^{k-1}a_ja_{k-j}
= \left(\sum_{j\geq 1} a_jq^j \right)^2 =1,
\end{equation}
while on the other hand, since $a_1=1$,
\begin{equation}
\sum_{k \geq 2} (k+1)a_k q^k =\sum_{k \geq 1} k a_k q^k +  
\sum_{k \geq 1} a_k q^k - 2q = g+1-2q.
\end{equation}
We obtain by this way $g=2q$, so that $x=1/2$, and thus $c_1=1/2$.

To conclude the proof, it suffices to check that 
$\sum_{k \geq 1} k^2 c_k=2$ with the previous values for $x$ and $g$.
But again, we obtain on the one hand that
\begin{equation}
\sum_{k\geq 2} k(k+1)a_k q^k =\sum_{k\geq 1} k^2 a_k q^k + 
\sum_{k\geq 1} k a_k q^k -2q = g \sum_{k\geq 1} k^2 c_k + g-2q,
\end{equation}
while on the other hand,
\begin{eqnarray}
\sum_{k\geq 2} k(k+1)a_k q^k =
\sum_{k\geq 2} k q^k \sum_{j=1}^{k-1}a_ja_{k-j}
= \sum_{j \geq 1} a_j q^j \sum_{k \geq j+1} k q^{k-j}a_{k-j} \ala
= \sum_{j \geq 1} a_j q^j \sum_{l \geq 1} (j+l) q^l a_l
=2 \left(\sum_{j \geq 1} a_j q^j\right)\left(\sum_{l \geq 1} l q^l a_l \right)
=2g.
\end{eqnarray}
Hence $g\sum_{k\geq 1} k^2 c_k +g - 2q = 2g$, so that 
$\sum_{k\geq 1} k^2 c_k = 1+2q/g=2$.
\end{proof}

To conclude this section, let us give some numerical results.

We obtain numerically, computing the values of $a_1,a_2,...,a_{10000}$,
and studying the function $z \mapsto \sum_{1}^{10000} a_k(z/2)^k$,
that $g \simeq 1.4458$, with quite a
good precision. 
We deduce then from Proposition \ref{sscf} that
at equilibrium, the mean-field model (\ref{coagfrag}) satisfies
\begin{eqnarray}
&c_1=0.5, \; c_2 \simeq 0.1204 , \;  c_3 \simeq 0.04354, \ala  
&c_4 \simeq 0.01679 ,\;  c_5 \simeq 0.006574, \;  c_6 \simeq 0.002582,\ala
&\sum_{k \geq 1} c_k \simeq 0.6916 ,\; \sum_{k \geq 1} k^2 c_k =2. \nonumber
\end{eqnarray}

On the other hand, 
simulating $10^8$ times the mass $M$ of the particle
containg the edge $(0,1)$ in the avalanche process $\eta_0$ at equilibrium
(see Appendix \ref{appendix}) we obtain the following 
Monte-Carlo approximations for $c_k(\eta_0):= P[M=k]/k$ 
\begin{eqnarray}
&c_1(\eta_{0}) \simeq 0.499934 , \; 
c_2(\eta_{0}) \simeq 0.12312 , \;
c_3(\eta_{0}) \simeq 0.0422142, \ala
&c_4(\eta_{0}) \simeq 0.0161849 ,  \;
c_5(\eta_{0}) \simeq 0.00648257 , \;
c_6(\eta_{0}) \simeq 0.00263739, \ala
&\sum_{k \geq 1} c_k(\eta_0) \simeq 0.692419 ,\; 
\sum_{k \geq 1} k^2 c_k(\eta_0) \simeq 1.99979. \nonumber
\end{eqnarray}

It appears clearly that the two sets of values are very similar,
even if numerical computations indicate
that no equality holds, except maybe concerning
$c_1$ and $\sum_{k} k^2 c_k$.
We have no explanation for this phenomenon.
It might indicate that correlations between the masses of 
successive clusters are nearly insignificant.

We have no proof that the mean-field model is the (very fast) 
limit, in some asymptotic regime, 
of the avalanche process.

\appendix
\section{Appendix}\label{appendix} \setcounter{equation}{0}

Let us now write down the simulation algorithm, which we deduce from
Sections \ref{reverse} and \ref{estimate}.
For $l\geq 0$, the algorithm below simulates a random variable 
$(\heta_0(i))_{i\in[-l,l]}$ with
law $\Pi_{[-l,l]}$, where $\Pi$ is the unique invariant distribution
of the avalanche process. 
The idea is to simulate $N$, $\zeta$ and $V$  in an a.s. finite random  
space-time domain and then to reconstruct $\eta$ following the graphical 
construction~\ref{gcab}. 

\vip

We construct a random process $\hzeta_n(k)$ 
containing the values of $\zeta$ at some random times $\hat T_n$ 
(times of jumps of $N$ and $V$ in a finite growing spatial domain
$[\ell_n,r_n]$) and an 
additionnal information: roughly, 

$\hzeta_{n}(k) = 0$ if $\zeta_{\hat T_n}(k)=0$;  

$\hzeta_{n}(k) = 2$ if $\zeta_{\hat T_n}(k)=1$ and 
$(k,\hat T_n)$ belongs to the box 
delimited by the contour processes; 

$\hzeta_{n}(k) = 1$ if $\zeta_{\hat T_n}(k)=1$ and $(k,\hat T_n)$ is outside
the box. 

\vip

We invite the reader to have a look to Figures 2 (for Step 1) 
and 4 (for Step 2) while reading
the simulation algorithm below. We will say {\it the box} 
for the box delimited by the contour processes.

\vip

{\bf Simulation Algorithm for $\Pi_{[-l,l]}$}

\vip

{\bf Step 0: Initialization.} 

Simulate the initial Bernoulli configuration $\zeta_0(k)$ 
for $k\in [l_0,r_0]$, where $l_0$ (resp. $r_0$) is the 
first vacant site on the left (resp. right) of $-l$ (resp. $l$).

If $\zeta_0(k)=0$ for all $k\in [-l,l]$, set $\heta_0(k)=0$ for all
$k\in [-l,l]$, and stop here.

Else, set $\hzeta_0(k)=2\zeta_0(k)$ for $k\in [l_0,r_0]$,
set $n=0$, and proceed to Step 1.

\vip

{\it Initially, all the sites in $[l_0+1, r_0-1]$ are in the box. We thus
assign the value $0$ to vacant sites, and the value $2$ to occupied sites}.

\vip

{\bf Step 1: Simulation of black/grey marks and contour processes} 

Set $n=n+1$. 

Choose $i_n$ uniformly in $[\ell_{n-1},r_{n-1}]$ 
{\it representing the involved site}.

Choose $m_n \sim$ Ber$(1/2)$, {\it here  
$m_n=0$ for a black mark, $m_n=1$ for a grey mark}.

\begin{itemize}
\item If $m_n = 0$, and $\hzeta_{n-1}(i_n)\geq 1$ then we  set 
$\hzeta_{n}(i_n) = 0$.

{\it If the site $i_n$ is occupied, it becomes vacant
due to a black mark.}

\item If $m_n=1$ and $\hzeta_{n-1}(i_n)=2$; 

set $\hzeta_{n}(i_n)=1$ if 
$\hzeta_{n-1}(i_n-1)=\hzeta_{n-1}(i_n+1)=0$ and  if
$\forall k\in [\ell_{n-1},i_n-1], \, \hzeta_{n-1}(k) \leq 1$ or 
$\forall k \in [i_n+1,r_{n-1}], \, \hzeta_{n-1}(k) \leq 1 $;

{\it The site $i_n$ remains occupied but 
leaves the box due to a grey mark, 
because its two neighbors are vacant, and because it is on the 
boundary of the box.}

otherwise, set $\hzeta_{n}(i_n)=2$.

{\it The site $i_n$ remains occupied and in the box,
because either the site $i_n$ is in the strict interior of the box
or one of its neighbors is occupied.}

\item If $m_n =0$ and  $\hzeta_{n-1}(i_n) = 0$ then we consider
the connected component $I_n$ of occupied sites (plus $i_n$) around $i_n$
(at time $n-1$).

If  
$\forall k \in I_n \cup [i_n+1, r_{n-1}], \, \hzeta_{n-1}(k) \leq 1$ or 
if $\forall k \in I_n \cup [\ell_{n-1}, i_n-1], \, \hzeta_n(k) \leq 1$, 
then, for all $k \in I_n$, set $\hzeta_{n}(i_n)=1$.

Otherwise, set  $\hzeta_{n+1}(k)=2$ for all $k \in I_n$.

{\it The site $i_n$ becomes occupied due to a black mark.
Then its whole connected component of occupied sites joins
the box, except if all these sites were outside the box at time $n-1$.}

\item Set $\hzeta_{n}(k) = \hzeta_{n-1}(k)$ for all sites 
$\ell_{n-1} \leq k \leq r_{n-1}$ of which the value (at time $n$) 
has not been defined yet.

{\it We update all other sites. Observe that we have not considered the case
$m_n=1,\hzeta_{n-1}(i_n)\in \{0,1\}$:
grey marks have no effect on vacant sites
in the Bernoulli process,
and cannot increase the number of sites in the box.}
\end{itemize}

\vip

If $\hzeta_n(r_{n-1})\leq 1$, set $r_n=r_{n-1}$. Else, 
consider $s^r_n\sim$ Geo$(1/2)$ (i.e.
$P[s^r_n=k]=(1/2)^k$ for $k\geq 1$), set $r_n=r_{n-1}+s^r_n$,
$\hzeta_{n}(k)=2$
for $k\in [r_{n-1},r_n-1]$, and $\hzeta_{n}(r_n)=0$.

{\it We extend the box if necessary, i.e. when
$\hzeta_n(r_{n-1})$ is occupied and in the box.
We thus extend the Bernoulli process to the right
until we meet a vacant site, at some site $r_n$. Then the
(occupied) sites $k\in [r_{n-1}+1,r_n-1]$ are in the box.}

\vip

If $\hzeta_n(\ell_{n-1})\leq 1$, set $\ell_n=\ell_{n-1}$. Else, 
consider  $s^\ell_n\sim$ Geo$(1/2)$, set $\ell_n=\ell_{n-1}-s^\ell_n$,
$\hzeta_{n}(k)=2$
for $k\in [\ell_n+1,\ell_{n-1}]$, and $\hzeta_{n}(\ell_n)=0$.

{\it Here we use the same arguments on the left of the domain.}

\vip

Check $\{ \ell_n \leq k \leq r_n, \hzeta_n(k) = 2\}$. 
If this set is non empty then repeat Step 1. Otherwise, set 
$T=n$ and proceed to Step 2. 

{\it If all the sites have a value equal to $0$ or $1$, this means
that the contour processes have met.}

\vip

{\bf Step 2: Deduction of the avalanche invariant realization.}

Start with $\heta_{T}(k)=0$ for all $\ell_T\leq k \leq r_T.$ 
Then for all $1\leq n \leq T$, 
define recursively (for  $n$ decreasing from $T$ to $1$) 
$\heta_{n-1}(k)$, for all $ \ell_{n-1}\leq k \leq r_{n-1}$, by 

\begin{itemize}
\item $\heta_{n-1}(k)=0$ if $m_n=0$ and if $k$ belongs to the 
connected component of occupied sites of $i_n$ (in $\heta_{n}$).

{\it Black marks kill connected components of occupied sites.}

\item $\heta_{n-1}(k)=1$ if $k=i_n$, 
if $m_n=1$ and if $\hzeta_{n-1}(i_n)\geq 1$.

{\it Grey marks make appear flocks when the Bernoulli process
is occupied.}

\item $\heta_{n-1}(k)=1$ if $k=i_n$, if $m_n=0$ and if 
$\heta_{n}(i_n)=\hzeta_{n}(i_n)=0$.

{\it Black marks make appear flocks at vacant sites.}

\item $\heta_{n-1}(k)=\heta_{n}(k)$ for all sites $k\in [\ell_{n-1},r_{n-1}]$ 
for which the value
$\heta_{n-1}(k)$ has not been defined yet.
\end{itemize}

\vip

{\bf Conclusion.}

Then $\{\heta_0(k), k\in [-l,l]\}$ is distributed 
according to $\Pi_{[-l,l]}$.

\vip

Remark that the number of iterations of Step 1 is finite
due to our results: when the contour processes $L^{-l}$ and
$R^l$ meet, there are no more sites with value $2$.

\vip

{\bf Alternative}

We finally propose another version of Step 1. The advantage 
is that the number of involved sites is much smaller. 
The idea is to take better 
advantage of the so-called {\it grey marks} : we will keep track 
only of what may really be needed to reconstruct $\{\heta_0(k), k\in [-l,l]\}$. 

\vip

{\bf Step 1'.}

Set $n=n+1$. Choose $i_n$ uniformly in $[\ell_{n-1},r_{n-1}]$.
Choose $m_n \sim$ Ber$(1/2)$,

\begin{itemize}
\item If $m_n = 0$, and $\hzeta_{n-1}(i_n)\geq 1$ then we  set 
$\hzeta_{n}(i_n) = 0$.

\item If $m_n=1$ and $\hzeta_{n-1}(i_n)=2$; we set
$\hzeta_{n}(i_n)=1$ as soon as $\hzeta_{n-1}(i_n-1)\leq 1$ or 
$\hzeta_{n-1}(i_n+1)\leq 1$; else we let $\hzeta_{n}(i_n)=2$.

\item If $m_n =0$ and  $\hzeta_{n-1}(i_n) = 0$ then we consider
\begin{eqnarray*}
I_n+=\{k\in [i_n+1,r_{n-1}], \hzeta_{n-1}(k)=2,\; \forall
i_n<j< k,  \hzeta_{n-1}(j)\geq 1\},\ala
I_n-=\{k\in [\ell_{n-1}, i_n-1], \hzeta_{n-1}(k)=2, \; \forall
k<j< i_n,  \hzeta_{n-1}(j)\geq 1\}.
\end{eqnarray*}
If $I_n+= I_n - = \emptyset$, we set $\hzeta_n(i_n)=1$.

Else, we set $\hzeta_n(k)=2$ for all $k \in [\min (I_n-),\max(I_n+)]$.

\item Set $\hzeta_{n}(k) = \hzeta_{n-1}(k)$ for all sites 
$\ell_{n-1} \leq k \leq r_{n-1}$ of which the value (at time $n$) 
has not been defined yet.
\end{itemize}

\vip

\vip

If there is $k \in [\ell_{n-1},r_{n-1}]$ such that $\hzeta_{n}(k)=2$
and $\hzeta_{n}(j)\geq 1$ for all $\ell_{n-1}\leq j \leq k$, 
consider $s_n^\ell \sim$ Geo$(1/2)$, set $\ell_n=\ell_{n-1}- s^\ell_n$,
$\hzeta_{n}(j)=1$ for $j \in [\ell_{n}+1, \ell_{n-1}-1]$ and $\hzeta_{n}(\ell_n)=0$.

Otherwise, set $\ell_n=\ell_{n-1}$. 

\vip 

Act symmetrically on the right
of the domain.

\vip

Check $\{ \ell_n \leq k \leq r_n, \hzeta_n(k) = 2\}$. 
If this set is non empty then repeat Step 1'. Otherwise, set 
$T=n$ and proceed to Step 2.

\vip

Let us emphasize the differences between Step 1 and Step 1'. 
The domain is extended only if there is $k$ with $\hzeta_n(k) = 2$ 
in the {\it connected} component touching the boundary. 
The first and fourth switching rules of $\hzeta_n(k)$ 
are left unchanged. In the second one, it is made {\it easier} to set 
$\hzeta_n(k)$ from $2$ to $1$. In the third one, a much smaller set 
of indices is switched from $1$ to $2$. 

\vip

One may understand that changing Step 1 into Step 1' does not
change the law of the final values $(\heta_0(k))_{k\in[-l,l]}$. 
It defines, in some sense, a much more precise contour process 
(leftmost and rightmost sites $k$ with $\hzeta_n(k) = 2$) than
the one defined in Section \ref{estimate}:
we numerically observe that using Step 1', the algorithm is
$15$ times faster than when using Step 1.
We hope to take advantage of
this idea to generalize our methods to more general particle systems.

\def\refname{References}

\end{document}